\documentclass[a4paper,11pt]{elsarticle}

\usepackage[T1]{fontenc}
\usepackage{newtxtext}
\usepackage{newtxmath}
\usepackage{amsmath}

\usepackage{amssymb}
\usepackage{bm}
\usepackage{graphicx}
\usepackage{float}
\usepackage{placeins}
\usepackage{booktabs}
\usepackage{enumitem}
\usepackage[margin=20mm]{geometry}
\usepackage{hyperref}
\usepackage{subcaption}
\usepackage{tikz}
\usetikzlibrary{arrows.meta}
\definecolor{SwanseaBlue}{RGB}{0,51,102}

\biboptions{sort&compress}
\graphicspath{{./figures/}}

\DeclareMathOperator{\tr}{tr}

\hypersetup{
    colorlinks=true,
    linkcolor=black,
    citecolor=black,
    urlcolor=blue
}
\begin{document}

\begin{frontmatter}

\title{On finite elements for geometrically-exact planar beams with hyperelastic material models}

\author[add1]{Abhishek Ghosh\corref{cor1}}
\ead{abhishek.ghosh@swansea.ac.uk}
\author[add2]{Chennakesava Kadapa\corref{cor1}}
\ead{c.kadapa@napier.ac.uk}
\author[add1]{Djordje Peri\'c}
\author[add1]{Mokarram Hossain\corref{cor1}}
\ead{mokarram.hossain@swansea.ac.uk}
\cortext[cor1]{Corresponding author}

\address[add1]{Zienkiewicz Institute for Modelling, Data and AI, Swansea University, Bay Campus, United Kingdom}
\address[add2]{School of Computing, Engineering and Built Environment,
Edinburgh Napier University, Edinburgh EH10 5DT, United Kingdom}

\begin{abstract}
A geometrically-exact planar beam finite-element framework is developed for the analysis of beams with hyperelastic constitutive models.
A conventional three-field beam model is first formulated with full geometrical nonlinearity and a linear elastic material law.
The kinematics are then enriched by an additional global field that permits deformation through the in-plane thickness.
The out-of-plane direction is incorporated separately through plane-strain and plane-stress constitutive reductions, allowing the cross-sectional deformation assumptions and constitutive behaviour to be formulated independently.
Numerical results for the present benchmarks show that the plane-stress condition gives essentially coincident centreline deformations for the Poisson's ratio values considered.
In contrast, the plane-strain condition, which is widely employed in the literature, makes the formulation sensitive to Poisson's ratio coupling, with the beam becoming significantly stiffer at higher Poisson's ratios. 
Comparisons between Saint Venant-Kirchhoff and Neo-Hookean hyperelastic models further demonstrate close agreement in the global deformation and in-plane thickness response, with small constitutive differences becoming more apparent in the local out-of-plane stretch.
The proposed approach offers a simple and consistent framework to integrate hyperelastic material models into planar geometrically-exact beam models.
\end{abstract}
\begin{keyword}
Geometrically-exact beam, Finite deformation, Hyperelasticity, Plane strain, Plane stress, Poisson coupling
\end{keyword}

\end{frontmatter}

\section{Introduction}
\label{sec:introduction}

Thin and slender structures occur throughout engineering, physics, and biology. Beams, rods, fibres, and filaments are characterised by a dominant longitudinal dimension, while plates and shells are thin relative to their in-plane dimensions. Such geometries are also central to soft and multifunctional systems, where large changes of shape often form part of the intended operation. Representative examples include compliant filaments undergoing coupled stretching, bending, shear, and torsion \cite{Gazzola2018}, soft robotic systems \cite{RusTolley2015}, dielectric and electro-active polymer devices \cite{Pelrine2000,DadgarRad2019,KadapaHossain2020Electro}, and magnetically responsive elastomers and slender hard-magnetic structures \cite{BastolaHossain2020,KadapaHossain2022MAP,YanAbbasiReis2022}. Magnetic actuation is particularly attractive because it permits remote, contactless control and has been exploited in reconfigurable soft devices \cite{KimYukZhao2018}, flexible endovascular probes \cite{Pancaldi2020}, and magnetically driven swimmers interacting with surrounding fluids \cite{RenHuDongSitti2019}. Related developments in large-strain electromechanical beam modelling further illustrate the increasing need for reduced structural descriptions that remain compatible with nonlinear material behaviour and coupled physical fields \cite{FirouziRabczukBonetZur2024}.

The separation of geometric scales in these problems naturally motivates dimensionally reduced structural models. A three-dimensional continuum resolves the deformation throughout the body, whereas beams or rods reduce the kinematics to a reference line supplemented by variables describing the cross-section \cite{SimoFoxRifai1989,CaliriFerreiraTita2016}. Here, ``reduced'' refers to this structural dimensional reduction and not to projection-based model-order reduction of an already discretised system \cite{BennerGugercinWillcox2015}. Three-dimensional models remain necessary when local through-thickness stresses, detailed cross-sectional warping, heterogeneous material distributions, or other genuinely three-dimensional features are of primary interest. For the global response of a highly slender structure, however, a beam model can represent the dominant axial, bending, and shear modes with substantially fewer unknowns, thereby enabling computationally efficient simulations. This economy becomes particularly valuable in coupled simulations, where structural deformation forms only one part of a larger multiphysics problem, as demonstrated, for example, in large-deformation fluid-structure interaction \cite{KadapaDettmerPeric2016}.

The foundations of finite-deformation beam and rod theory precede modern geometrically exact finite-element formulations. Early developments in one-dimensional nonlinear elasticity include the hyperelastic rod theory of Antman and Warner \cite{AntmanWarner1966} and Reissner's finite-strain theory for the planar beam problem, in which axial extension, transverse shear, and finite rotation are treated within a consistent nonlinear framework \cite{Reissner1972}. The formulation of Simo \cite{Simo1985} and the computational developments of Simo and Vu-Quoc \cite{SimoVuQuoc1986} subsequently provided a systematic treatment of finite-rotation kinematics, variational equilibrium, and consistent linearisation. For planar beams, Saje developed variational and finite-element formulations for straight and curved members undergoing finite deformation \cite{Saje1990,Saje1991}, while Ibrahimbegovi\'c and Frey examined linear and geometrically nonlinear planar beam elements together with membrane- and shear-locking behaviour \cite{IbrahimbegovicFrey1993}. Complementary finite-element developments addressed nonlinear beam dynamics with finite rotations \cite{CardonaGeradin1988} and geometrically nonlinear Reissner-type curved beam elements \cite{Ibrahimbegovic1995}. Simo and Vu-Quoc further extended geometrically exact rod kinematics to include shear and torsion-warping deformation \cite{SimoVuQuoc1991}. Later developments examined objectivity and frame indifference \cite{JelenicCrisfield1999,BetschSteinmann2002,RomeroArmero2002}, addressed the large-deformation analysis of flexible beams \cite{PaiPalazotto1996}, clarified connections with three-dimensional finite elasticity \cite{AuricchioCarotenutoReali2008}, and developed strain-based finite-element formulations of geometrically exact beams \cite{ZupanSaje2003}, alongside robust displacement--rotation formulations for slender structures undergoing arbitrarily large changes of configuration \cite{MeierPoppWall2014}.

It is, however, essential to distinguish between geometric and material nonlinearities. A beam may exhibit large displacements and finite rotations while still obeying a linear material law between its generalised strains and section forces/moments. In contrast, soft solids often undergo substantial material strains, requiring their constitutive behaviour to be characterised by hyperelastic strain-energy functions. To integrate such a material description into a reduced beam model, one must ensure compatibility between the deformation captured by the beam kinematics and the three-dimensional strain measures used in the constitutive equations. At the continuum level, large-strain formulations based on polyconvex strain-energy functions offer a systematic framework for nonlinear hyperelasticity \cite{BonetGilOrtigosa2015}. For rubber-like materials, a broad range of phenomenological and micromechanically motivated hyperelastic models have been developed together with the associated stress measures and consistent material tangent operators required for finite-element implementation \cite{SteinmannHossainPossart2012,HossainSteinmann2013}. Network-based constitutive descriptions, including eight-chain and full-network models, provide a complementary micromechanical representation of rubber hyperelasticity \cite{HossainAminKabir2015}. A corresponding bridge between three-dimensional polyconvex elasticity and geometrically exact beam theories was constructed in \cite{OrtigosaGilBonetHesch2016}, where a direct relationship between three-dimensional elasticity and classical beam strain measures was derived. Related constitutive reductions have also been developed within special Cosserat rod theory. The Helical Cauchy--Born construction provides a route for obtaining nonlinear rod constitutive relations from the underlying three-dimensional elasticity \cite{AjeetKumarKumarGupta2016}, while a computational extension of this approach permits arbitrary three-dimensional hyperelastic material behaviour and determines the resulting force, moment, and stiffness response through a nonlinear cross-sectional warping problem \cite{AroraAjeetKumarSteinmann2019}. Alternative strategies include beam theories with deformable or extensible directors \cite{ChoiSauerKlinkel2021}, formulations incorporating explicit through-thickness variables \cite{ZurFirouziRabczukZhuang2023}, and hyperelastic models tailored to stimulus-responsive materials \cite{DadgarRad2019,DadgarRadMahjoubHossain2024,FirouziRabczukBonetZur2024}.

The treatment of the cross-section becomes especially important for soft polymeric materials, which are commonly modelled as incompressible or nearly incompressible. This issue is well established in electro-active and magneto-active continua \cite{KadapaHossain2020Electro,BastolaHossain2020,KadapaHossain2022MAP}, while mixed displacement-pressure formulations provide one route for treating nearly and fully incompressible three-dimensional hyperelasticity \cite{KadapaHossain2022Mixed}. Dimensionally-reduced beam models introduce additional difficulties because deformation components suppressed by their kinematic assumptions can remain constitutively coupled to the retained components. The mechanical significance of normal cross-sectional extension was recognised by Naghdi and Rubin \cite{NaghdiRubin1989}, who showed that excluding transverse normal deformation may remove relevant features of the structural response. More recent hyperelastic beam formulations have therefore introduced additional cross-sectional degrees of freedom or thickness variables \cite{HeLouDongKitipornchaiYang2018PlaneStrain,HeLouDongKitipornchaiYang2018Shearable,DadgarRad2019,ChoiSauerKlinkel2021,ZurFirouziRabczukZhuang2023}. A related reduction for strip-like structures has represented the cross-sectional deformation through variables describing the deformed cross-sectional centreline and the orientation of the thickness lines, allowing nonlinear constitutive response to be retained within a one-dimensional special Cosserat rod description \cite{RoushanKumarAgarwalAjeetKumar2024}.

The presence of a transverse kinematic field does not, by itself, determine how a reduced model accommodates three-dimensional constitutive coupling. For a linear isotropic material, Young's modulus and Poisson's ratio determine the Lam\'e parameters. In the special case of zero Poisson's ratio, with or without the additional thickness stress unknown, the associated coupling between normal and shear strains disappears. Ortigosa et al.\ \cite{OrtigosaGilBonetHesch2016} discussed the relation between beam constitutive restrictions and this special case, while Choi et al.\ \cite{ChoiSauerKlinkel2021} introduced extensible directors together with enhanced strains to alleviate Poisson locking when general three-dimensional constitutive laws are employed. Thickness-stretchable finite-strain beams have also been developed explicitly under plane-strain assumptions \cite{HeLouDongKitipornchaiYang2018PlaneStrain,HeLouDongKitipornchaiYang2018Shearable}. These contributions show that cross-sectional kinematics and the constitutive restriction in the remaining transverse direction cannot, in general, be considered independently.

When planar beam models are used with three-dimensional hyperelastic laws, the treatment of strain and stress constraints in the plane-strain and plane-stress models becomes particularly important. Under the plane-strain model, as the deformation in the out-of-plane direction is suppressed, an out-of-plane stress will develop as the reaction required to enforce that kinematic constraint. Under the plane-stress model, the out-of-plane stress is assumed to vanish and the corresponding out-of-plane deformation must instead be determined consistently from the constitutive response. Klinkel and Govindjee \cite{KlinkelGovindjee2002} developed a general local integration-point procedure for incorporating three-dimensional finite-strain constitutive models into beam and shell elements subject to prescribed zero-stress conditions, together with the corresponding consistent condensed material tangent. In contrast, the present formulation introduces an independent global field to represent deformation through the in-plane beam thickness, while the remaining out-of-plane constitutive response is treated separately through either a plane-strain or plane-stress reduction. In the plane-stress case, the required out-of-plane stretch is determined locally from the zero-stress condition and condensed at material level, without introducing an additional global beam degree of freedom. Related local plane-stress reductions are well established for general hyperelastic shells \cite{KiendlHsuWuReali2015} and have been formulated systematically for arbitrary compressible hyperelastic materials \cite{AhmadiMcBrideSteinmannSaxena2025}. Fully incompressible hyperelastic beam formulations have also combined transverse plane-stress conditions with an explicit incompressibility constraint \cite{DadgarRadSahraee2021}. Plane stress and incompressibility are therefore distinct restrictions: the former permits constitutive relaxation in the out-of-plane direction, whereas exact incompressibility additionally requires preservation of local volume.

In this context, the goal of the present work is to integrate in-plane thickness deformation, Poisson-dependent constitutive coupling, and the distinction between plane-strain and plane-stress reductions within a common geometrically exact planar beam framework. This allows the influence of Poisson's ratio to be examined without restricting the constitutive description to the special zero-Poisson-ratio case, \emph{which is the classical case predominantly used in the literature on finite elements for planar beam problems but not practically meaningful for hyperelastic cases}. The comparison is performed using identical enriched beam kinematics, so that differences between the plane-strain and plane-stress responses arise from the treatment of the out-of-plane constitutive behaviour. Particular attention is devoted to nearly incompressible materials, for which restricting transverse deformation can strongly influence the response of the planar beam models.

The paper is organised as follows. Section~\ref{sec:linear_beam} presents the standard three-field model and assesses it against a three-dimensional continuum solution under a vertical tip load and the analytical pure-bending solution under an applied end moment. Section~\ref{sec:hyperelastic_thickness} develops the enriched kinematics, finite-element treatment, plane-strain and plane-stress reductions, and the corresponding end-moment comparisons. The Saint-Venant-Kirchhoff and Neo-Hookean material specialisations used in these benchmarks are collected in \ref{app:material_specialisations}. Concluding remarks are provided in Section~\ref{sec:conclusions}.

\section{Geometrically-exact planar beam formulation with linear elastic material}
\label{sec:linear_beam}
The formulation developed in this section is geometrically nonlinear but materially linear. Large displacements and finite rotations are retained through the finite-deformation kinematics, so that the generalized beam
strains depend nonlinearly on the centreline displacements and cross-sectional rotation, while the section constitutive law remains linear in these strain measures. Within the adopted conventional beam description, no independent transverse deformation of the cross-section is introduced. The present section adopts a Hooke-type linear elastic material response. In Section~\ref{sec:hyperelastic_thickness}, the kinematics are enriched by an independent thickness-stretch field and coupled to general hyperelastic constitutive laws, so that transverse deformation and material nonlinearity can be treated within the same framework.

\subsection{Formulation} \label{sec:linear_formulation}
The formulation is developed systematically from the beam geometry and kinematics through the strain measures, constitutive relations, variational statement, and finite-element discretisation.

\subsubsection{Beam geometry and kinematics}
\begin{figure}[H]
\centering
\begin{tikzpicture}[>=Stealth, scale=1.05]

\draw[->] (0,0) -- (10,0) node[below right] {$X$};
\draw[->] (0,0) -- (0,7) node[above left] {$Y$};

\draw[thick] (1.5, 1.0) -- (7, 1.0);
\draw[thick] (1.5, 1.8) -- (7, 1.8);
\draw[thick] (7, 1.0) -- (7, 1.8);

\draw[thick]
(1.5, 1.0)
.. controls (1.3, 1.2) and (1.7, 1.6) ..
(1.5, 1.8);

\draw[dash dot] (1.0, 1.4) -- (7.5, 1.4);

\node at (4.25, 1.4)
[below, font=\small]
{Original configuration};

\draw[thick]
(1.5, 3.5)
to[out=5, in=195]
(7.8, 4.5);

\draw[thick]
(1.5, 4.5)
to[out=5, in=195]
(8.1, 5.5);

\draw[thick] (7.8, 4.5) -- (8.1, 5.5);

\draw[thick]
(1.5, 3.5)
.. controls (1.3, 3.7) and (1.7, 4.3) ..
(1.5, 4.5);

\draw[dash dot]
(1.0, 4.0)
to[out=5, in=195]
(8.5, 5.15);

\node at (4.5, 4.32)
[rotate=9, below, font=\small]
{Deformed configuration};

\draw[thin] (7, 1.8) -- (7, 3.5);
\draw[thin] (7.95, 3.5) -- (7.95, 4.5);

\draw[<->]
(7, 3.8) -- (7.95, 3.8)
node[midway, above] {$u$};

\draw[thin] (7.95, 5.0) -- (9.5, 5.0);

\draw[<->]
(9.0, 0) -- (9.0, 5.0)
node[midway, right] {$v$};

\draw[thin]
(7.95, 5.0) -- (7.95, 7.0);

\draw[thin]
(7.95, 5.0) -- ++(73.3:2.2);

\draw[<->]
(7.95, 5.0) +(90:1.8)
arc (90:73.3:1.8)
node[midway, above, sloped, yshift=-1pt] {$\beta$};

\end{tikzpicture}

\caption{Initial and displaced configurations of the beam for a three-field formulation.}
\label{fig:beam_config}
\end{figure}
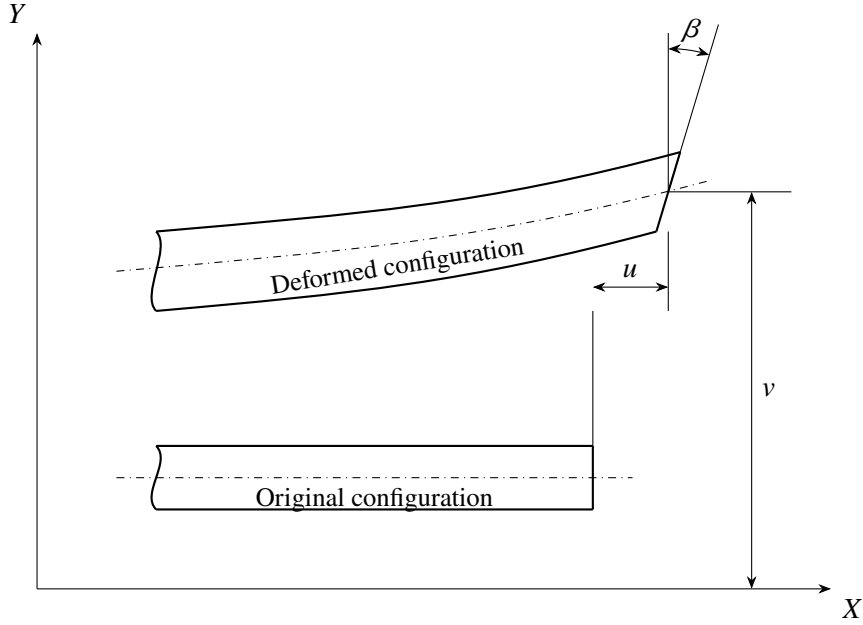

For a planar beam occupying the reference domain
$\Omega_0=[0,L]\times A$, we choose $X\in[0,L]$ as the material coordinate
along the centreline and $Y$ as the coordinate across the cross-section. With $u$ and $v$ as the centreline translations in the $X$- and $Y$-directions and $\beta$ as the cross-sectional rotation, as illustrated in Figure~\ref{fig:beam_config}, the current coordinates of a material point initially located at $(X,Y)$ are
\begin{subequations}
\label{eq:coordinates_linear}
\begin{align}
x(X,Y)
&=
X+u(X)+Y\sin\beta(X),
\label{eq:x_linear}
\\
y(X,Y)
&=
v(X)+Y\cos\beta(X).
\label{eq:y_linear}
\end{align}
\end{subequations}

Equivalently, the displacement map can be written as
\begin{subequations}
\label{eq:placement_map_linear}
\begin{equation}
\boldsymbol{\varphi}(X,Y)
=
\boldsymbol{\varphi}_c(X)+Y\mathbf d(X),
\label{eq:placement_linear}
\end{equation}
with
\begin{equation}
\boldsymbol{\varphi}_c(X)
=
\begin{bmatrix}
X+u(X)\\
v(X)
\end{bmatrix},
\qquad
\mathbf d(X)
=
\begin{bmatrix}
\sin\beta(X)\\
\cos\beta(X)
\end{bmatrix}.
\label{eq:centerline_director}
\end{equation}
\end{subequations}
Here, $\mathbf d(X)$ is the unit director.

We can now compute the in-plane deformation gradient as
\begin{equation}
\mathbf F_{2D}
=
\begin{bmatrix}
F_{11} & F_{12}\\
F_{21} & F_{22}
\end{bmatrix}
=
\begin{bmatrix}
1+u_{,X}+Y\cos\beta\,\beta_{,X} & \sin\beta\\
v_{,X}-Y\sin\beta\,\beta_{,X} & \cos\beta
\end{bmatrix}.
\label{eq:F_linear}
\end{equation}
Here, a comma followed by a coordinate denotes differentiation with respect to that reference coordinate. For example, $u_{,X}=\dfrac{\partial u}{\partial X}$.
In the displacement Equation~\eqref{eq:placement_linear}, the transverse coordinate $Y$ is carried by the unit director $\mathbf d(X)$, so that the cross-section may translate and rotate but does not undergo an independent normal deformation in the $Y$-direction. 

The determinant of the deformation gradient becomes
\begin{subequations}
\label{eq:jacobian_linear}
\begin{equation}
J_{2D}
=
\det\mathbf F_{2D}
=
\Lambda+Y\beta_{,X}.
\label{eq:J_linear}
\end{equation}
Here,
\begin{equation}
\Lambda
=
(1+u_{,X})\cos\beta-v_{,X}\sin\beta.
\label{eq:Lambda_linear}
\end{equation}
\end{subequations}

\subsubsection{Generalized strain measures}
The beam strains follow from the Green-Lagrange strain tensor
\begin{equation}
\mathbf E
=
\frac12\left(\mathbf F^{\mathrm T}\mathbf F-\mathbf I\right).
\label{eq:green_lagrange_linear}
\end{equation}
Within the standard first-order through-thickness beam approximation, the
axial strain is retained up to terms linear in $Y$. The non-zero in-plane strain
components used by the beam model are therefore
\begin{equation}
E_{11}=E_0+YK_b,
\qquad
2E_{12}=\Gamma,
\qquad
E_{22}=0,
\label{eq:strain_distribution_linear}
\end{equation}
where
\begin{subequations}
\label{eq:generalized_measures_linear}
\begin{align}
E_0
&=
u_{,X}
+
\frac12\left(u_{,X}^2+v_{,X}^2\right),
\label{eq:E0_linear}
\\
\Gamma
&=
(1+u_{,X})\sin\beta
+
v_{,X}\cos\beta,
\label{eq:Gamma_linear}
\\
K_b
&=
\beta_{,X}\Lambda.
\label{eq:Kb_linear}
\end{align}
\end{subequations}
Here, $E_0$ is the centreline axial strain, $\Gamma$ is the transverse shear
strain measure, and $K_b$ is the bending curvature measure.

The generalized beam strain vector is defined as
\begin{equation}
\boldsymbol{\eta}
=
\begin{bmatrix}
E_0\\
\Gamma\\
K_b
\end{bmatrix}.
\label{eq:eta_linear}
\end{equation}

\noindent
\textbf{Remark}: The model assumes a plane-strain condition by default \cite{Zienkiewicz-book}, resulting in $F_{33}=1$, $E_{33}=0$, and $J=J_{2D}$. This restriction represents a boundary-value problem in which out-of-plane deformation is constrained; it is not appropriate when no such constraint is present. Together with the absence of independent in-plane normal deformation, it limits the representation of transverse deformation and the associated constitutive coupling in hyperelastic beams that are free to deform transversely.

\subsubsection{Strain generators and constitutive quantities}

A convenient tensor representation of the generalized strains is obtained by
introducing fixed strain generators. Let $\mathbf a_1$ and $\mathbf a_2$ denote
the Cartesian basis vectors associated with the reference $X$- and
$Y$-directions. The retained Green-Lagrange strain tensor can be written as
\begin{subequations}
\label{eq:strain_generators_linear}
\begin{equation}
\mathbf E(X,Y)
=
E_0\mathbf G_{E_0}
+
\Gamma\mathbf G_{\Gamma}
+
K_b\mathbf G_{K_b},
\label{eq:strain_generator_expansion}
\end{equation}
with
\begin{align}
\mathbf G_{E_0}
&=
\mathbf a_1\otimes\mathbf a_1,
\label{eq:G_E0}
\\
\mathbf G_{\Gamma}
&=
\frac12\left(
\mathbf a_1\otimes\mathbf a_2
+
\mathbf a_2\otimes\mathbf a_1
\right),
\label{eq:G_Gamma}
\\
\mathbf G_{K_b}
&=
Y\,\mathbf a_1\otimes\mathbf a_1.
\label{eq:G_Kb}
\end{align}
\end{subequations}
Thus, $E_0$ acts on the axial component, $\Gamma$ produces the symmetric
shear components $E_{12}=E_{21}=\Gamma/2$, and $K_b$ gives the linear
through-thickness variation of the axial strain.

To derive the section quantities variationally, let $\Psi(\mathbf E)$ denote
the strain-energy density and define the work-conjugate second
Piola-Kirchhoff stress and material tangent as
\begin{equation}
\mathbf S
=
\frac{\partial\Psi}{\partial\mathbf E},
\qquad
\mathbb L
=
\frac{\partial\mathbf S}{\partial\mathbf E},
\qquad
L_{IJKL}
=
\frac{\partial S_{IJ}}{\partial E_{KL}}.
\label{eq:constitutive_general_linear_section}
\end{equation}
These relations provide the section-level variational structure. For the
present three-field beam they are specialised below to the linear elastic law;
a general hyperelastic constitutive treatment is deferred to
Section~\ref{sec:hyperelastic_thickness}, where thickness stretch is included.

\subsubsection{Internal energy and section resultants}

The internal energy is
\begin{subequations}
\label{eq:energy_section_definition}
\begin{equation}
\Pi_{\mathrm{int}}(\mathbf q)
=
\int_{\Omega_0}\Psi(\mathbf E)\,\mathrm dV_0
=
\int_0^L\Phi\!\left(X;\boldsymbol{\eta}(\mathbf q)\right)\,\mathrm dX,
\label{eq:Pi_int_general_section}
\end{equation}
where
\begin{equation}
\Phi(X;\boldsymbol{\eta})
=
\int_A\Psi\bigl(\mathbf E(X,Y)\bigr)\,\mathrm dA
\label{eq:Phi_section}
\end{equation}
\end{subequations}
is the strain energy per unit reference length.

Using the strain-generator representation, the first variation becomes
\begin{align}
\delta\Pi_{\mathrm{int}}
&=
\int_{\Omega_0}\delta\mathbf E:\mathbf S\,\mathrm dV_0
\nonumber\\
&=
\int_0^L
\left(
\delta E_0\,T^p
+
\delta\Gamma\,S^p
+
\delta K_b\,M^b
\right)\mathrm dX.
\label{eq:delta_Pi_section}
\end{align}
The axial force, shear force, and bending moment resultants are therefore
\begin{subequations}
\label{eq:section_resultants_general}
\begin{align}
T^p
&=
\int_A\mathbf G_{E_0}:\mathbf S\,\mathrm dA
=
\int_A S_{11}\,\mathrm dA,
\label{eq:Tp_general}
\\
S^p
&=
\int_A\mathbf G_{\Gamma}:\mathbf S\,\mathrm dA
=
\int_A S_{12}\,\mathrm dA,
\label{eq:Sp_general}
\\
M^b
&=
\int_A\mathbf G_{K_b}:\mathbf S\,\mathrm dA
=
\int_A YS_{11}\,\mathrm dA.
\label{eq:Mb_general}
\end{align}
\end{subequations}
They are collected as
\begin{equation}
\mathbf r
=
\frac{\partial\Phi}{\partial\boldsymbol{\eta}}
=
\begin{bmatrix}
T^p\\
S^p\\
M^b
\end{bmatrix}.
\label{eq:r_section}
\end{equation}

\subsubsection{Consistent section tangent}

The consistent section tangent follows by linearising the resultants with
respect to the generalized beam strains. For an increment
$D\boldsymbol{\eta}$,
\begin{subequations}
\label{eq:consistent_section_tangent}
\begin{equation}
D\mathbf E
=
D\eta_{\alpha}\,\mathbf G_{\alpha},
\qquad
D\mathbf S
=
\mathbb L:D\mathbf E,
\label{eq:DE_DS_section}
\end{equation}
where the index $\alpha$ runs over the three generalized strain modes. Hence,
\begin{equation}
D r_{\alpha}
=
H_{\alpha\beta}D\eta_{\beta},
\qquad
H_{\alpha\beta}(X)
=
\int_A
\mathbf G_{\alpha}:\mathbb L:\mathbf G_{\beta}\,\mathrm dA.
\label{eq:H_general_section}
\end{equation}
Using the minor and major symmetries of $\mathbb L$, the section tangent can
be written as
\begin{equation}
\mathbf H
=
\int_A
\begin{bmatrix}
L_{1111} & L_{1112} & YL_{1111}\\[3pt]
L_{1112} & L_{1212} & YL_{1211}\\[3pt]
YL_{1111} & YL_{1211} & Y^2L_{1111}
\end{bmatrix}
\,\mathrm dA.
\label{eq:H_components_general}
\end{equation}
\end{subequations}
For a centroidal section and a material response that is uniform over the
cross-section, terms containing odd powers of $Y$ vanish.

\subsubsection{Linear elastic material}

The section quantities are now specialised to the materially linear law used with the conventional three-field beam. The retained axial and shear
stress components are taken as
\begin{subequations}
\label{eq:linear_material_relations}
\begin{equation}
S_{11}=E_YE_{11},
\qquad
S_{12}=2\mu E_{12}=\mu\Gamma,
\label{eq:linear_stress_law}
\end{equation}
where $E_Y$ is Young's modulus, $\nu$ is Poisson's ratio, and
\begin{equation}
\mu
=
\frac{E_Y}{2(1+\nu)}
\label{eq:G_linear}
\end{equation}
\end{subequations}
is the shear modulus.

For a centroidal reference cross-section, define
\begin{equation}
A=\int_A \mathrm dA,
\qquad
\int_A Y\,\mathrm dA=0,
\qquad
I=\int_A Y^2\,\mathrm dA.
\label{eq:section_properties}
\end{equation}
Using $E_{11}=E_0+YK_b$ and the centroidal condition in
Equation~\eqref{eq:section_properties} gives
\begin{equation}
\begin{aligned}
T^p
&=
E_YA\,E_0,
\\
S^p
&=
\mu A\,\Gamma,
\\
M^b
&=
E_YI\,K_b.
\end{aligned}
\end{equation}
A shear correction factor $\kappa_s$ is retained at section level to
illustrate how it enters a materially linear but geometrically nonlinear beam
formulation. Throughout this study, $\kappa_s=1$. The shear resultant is therefore
\begin{equation}
S^p
=
\kappa_s \mu A\,\Gamma.
\label{eq:shear_corrected_resultant}
\end{equation}
The linear elastic section law is
\begin{subequations}
\label{eq:linear_section_constitutive_set}
\begin{equation}
\mathbf r
=
\mathbf H_{\mathrm{lin}}\boldsymbol{\eta},
\label{eq:linear_section_law}
\end{equation}
with
\begin{equation}
\mathbf H_{\mathrm{lin}}
=
\begin{bmatrix}
E_YA & 0 & 0\\
0 & \kappa_s \mu A & 0\\
0 & 0 & E_YI
\end{bmatrix}.
\label{eq:H_linear}
\end{equation}
The corresponding strain energy per unit reference length is
\begin{equation}
\Phi_{\mathrm{lin}}
=
\frac12\boldsymbol{\eta}^{\mathrm T}
\mathbf H_{\mathrm{lin}}
\boldsymbol{\eta}
=
\frac12\left(
E_YA\,E_0^2
+
\kappa_s \mu A\,\Gamma^2
+
E_YI\,K_b^2
\right).
\label{eq:Phi_linear}
\end{equation}
\end{subequations}
The section law is linear in $\boldsymbol{\eta}$, but the beam problem remains geometrically nonlinear because $E_0$, $\Gamma$, and $K_b$ are nonlinear functions of $u$, $v$, and $\beta$.

\subsubsection{Weak form}

The equilibrium equations are obtained from the principle of virtual work.
For configuration-independent distributed line loads, define
\begin{subequations}
\label{eq:external_loading_linear}
\begin{equation}
\overline{\mathbf b}(X)
=
\begin{bmatrix}
\overline b_x(X)\\
\overline b_y(X)\\
\overline m(X)
\end{bmatrix},
\label{eq:distributed_load_linear}
\end{equation}
and let $\overline{\mathbf t}$ denote prescribed force and moment resultants
on the Neumann boundary $\Gamma_N$. The external virtual work is
\begin{equation}
\delta\Pi_{\mathrm{ext}}
=
\int_0^L
\delta\mathbf q^{\mathrm T}\overline{\mathbf b}\,\mathrm dX
+
\left.
\delta\mathbf q^{\mathrm T}\overline{\mathbf t}
\right|_{\Gamma_N}.
\label{eq:delta_Pi_ext_linear}
\end{equation}
\end{subequations}
The weak form is therefore
\begin{equation}
\delta\Pi
=
\delta\Pi_{\mathrm{int}}-\delta\Pi_{\mathrm{ext}}
=
0
\label{eq:weak_form_linear}
\end{equation}
for every admissible virtual field $\delta\mathbf q$.

\subsubsection{Finite-element discretisation and numerical implementation}

The finite-element equations follow by interpolating the three mechanical
fields over an element $\Omega_e=[0,L_e]$. With $n_e$ element nodes,
\begin{subequations}
\label{eq:FE_kinematics_linear}
\begin{equation}
\mathbf q^h(X)
=
\sum_{a=1}^{n_e}N_a(X)\mathbf q_a,
\qquad
\mathbf q_a
=
\begin{bmatrix}
u_a & v_a & \beta_a
\end{bmatrix}^{\mathrm T},
\label{eq:FE_interpolation_linear}
\end{equation}
and
\begin{equation}
\mathbf q_e
=
\begin{bmatrix}
u_1 & v_1 & \beta_1 & \cdots & u_{n_e} & v_{n_e} & \beta_{n_e}
\end{bmatrix}^{\mathrm T}.
\label{eq:q_element_linear}
\end{equation}
\end{subequations}

The variations of the generalized strains are
\begin{subequations}
\label{eq:strain_variations_linear}
\begin{align}
\delta E_0
&=
(1+u_{,X})\,\delta u_{,X}
+
v_{,X}\,\delta v_{,X},
\label{eq:delta_E0_linear}
\\
\delta\Gamma
&=
\sin\beta\,\delta u_{,X}
+
\cos\beta\,\delta v_{,X}
+
\Lambda\,\delta\beta,
\label{eq:delta_Gamma_linear}
\\
\delta K_b
&=
\Lambda\,\delta\beta_{,X}
-
\beta_{,X}\Gamma\,\delta\beta
+
\beta_{,X}
\left(
\cos\beta\,\delta u_{,X}
-
\sin\beta\,\delta v_{,X}
\right).
\label{eq:delta_Kb_linear}
\end{align}
\end{subequations}
Thus,
\begin{subequations}
\label{eq:B_matrix_linear_set}
\begin{equation}
\delta\boldsymbol{\eta}
=
\mathbf B\,\delta\mathbf q_e,
\qquad
\mathbf B
=
\begin{bmatrix}
\mathbf B_1 & \mathbf B_2 & \cdots & \mathbf B_{n_e}
\end{bmatrix},
\label{eq:B_linear}
\end{equation}
with nodal block
\begin{equation}
\mathbf B_a
=
\begin{bmatrix}
(1+u_{,X})N_{a,X} & v_{,X}N_{a,X} & 0\\[4pt]
\sin\beta\,N_{a,X} & \cos\beta\,N_{a,X} & \Lambda N_a\\[4pt]
\beta_{,X}\cos\beta\,N_{a,X} & -\beta_{,X}\sin\beta\,N_{a,X} &
\Lambda N_{a,X}-\Gamma\beta_{,X}N_a
\end{bmatrix}.
\label{eq:B_node_linear}
\end{equation}
\end{subequations}

The element internal force vector is
\begin{subequations}
\label{eq:element_residual_construction_linear}
\begin{equation}
\mathbf F_{\mathrm{int},e}
=
\int_0^{L_e}\mathbf B^{\mathrm T}\mathbf r\,\mathrm dX.
\label{eq:Fint_linear}
\end{equation}
For distributed loading, introduce the interpolation matrix
\begin{equation}
\mathbf N_q
=
\begin{bmatrix}
N_1 & 0 & 0 & \cdots & N_{n_e} & 0 & 0\\
0 & N_1 & 0 & \cdots & 0 & N_{n_e} & 0\\
0 & 0 & N_1 & \cdots & 0 & 0 & N_{n_e}
\end{bmatrix},
\label{eq:Nq_linear}
\end{equation}
so that
\begin{equation}
\mathbf F_{\mathrm{ext},e}^{\,b}
=
\int_0^{L_e}
\mathbf N_q^{\mathrm T}\overline{\mathbf b}\,\mathrm dX.
\label{eq:Fext_distributed_linear}
\end{equation}
Equivalent nodal contributions from concentrated forces and moments are added
to obtain $\mathbf F_{\mathrm{ext},e}$.

The element residual is
\begin{equation}
\mathbf{\mathcal{R}}_e
=
\mathbf F_{\mathrm{int},e}-\mathbf F_{\mathrm{ext},e}.
\label{eq:residual_linear}
\end{equation}
\end{subequations}
Its consistent linearisation gives
\begin{subequations}
\label{eq:tangent_decomposition_linear}
\begin{equation}
\mathbf{K}_e
=
\mathbf K_{\mathrm M}+\mathbf K_{\mathrm G},
\label{eq:Ke_linear}
\end{equation}
where
\begin{equation}
\mathbf{K}_{\mathrm M}
=
\int_0^{L_e}
\mathbf B^{\mathrm T}\mathbf H_{\mathrm{lin}}\mathbf B\,\mathrm dX
\label{eq:KM_linear}
\end{equation}
and
\begin{equation}
\mathbf{K}_{\mathrm G}
=
\int_0^{L_e}
\left(D\mathbf B^{\mathrm T}\right)\mathbf r\,\mathrm dX.
\label{eq:KG_linear}
\end{equation}
\end{subequations}
Following assembly, the Newton correction is obtained from
\begin{equation}
\mathbf K\,\Delta\mathbf q
=
-\mathbf{{\mathcal{R}}}.
\label{eq:Newton_linear}
\end{equation}

The element stiffness matrix and residual vectors are evaluated by numerical integration along the beam axis. Reduced axial integration is employed to alleviate shear locking in the slender-beam limit \cite{Saje1991,IbrahimbegovicFrey1993,ZienkiewiczTaylorToo1971}.


\subsection{Numerical examples}
\label{sec:linear_numerical_examples}

All beam simulations use incremental load control with an initial
load-factor increment of $0.05$, corresponding to a nominal subdivision
of the total load into 20 equal increments. At each attempted load step,
the equilibrium equations are solved using Newton-Raphson iterations,
with a maximum of 10 iterations. Convergence is accepted when the
absolute Euclidean norm of the assembled equilibrium residual restricted
to the free degrees of freedom falls below $10^{-8}$.
If this criterion is not satisfied within the iteration limit, the load
increment is halved for the subsequent attempt. For all models
considered, quadratic convergence is observed near equilibrium.

The conventional three-field formulation is assessed for a cantilever beam 
subjected to two loading cases: a vertical force at the free end and a
prescribed end moment. The length of the beam is $L=100\,\mathrm{mm}$, the in-plane height is $H=2\,\mathrm{mm}$, and the out-of-plane width is $W=1\,\mathrm{mm}$. 
The corresponding cross-sectional quantities are
\begin{equation}
A=WH=2\,\mathrm{mm}^2,
\qquad
I=\frac{WH^3}{12}=0.6666666667\,\mathrm{mm}^4.
\label{eq:linear_benchmark_section}
\end{equation}
The Young's modulus is $E_Y=0.606\,\mathrm{N/mm}^2=6.06\times10^5\,\mathrm{Pa}$. The centreline is discretised using 20 equal two-node elements, giving an element length of $5\,\mathrm{mm}$ and 21 nodes. The calculations are performed for different Poisson's-ratio values. The corresponding global centrelines are essentially coincident.

\subsubsection{Cantilever beam with vertical end load}
\label{sec:linear_end_point_load}

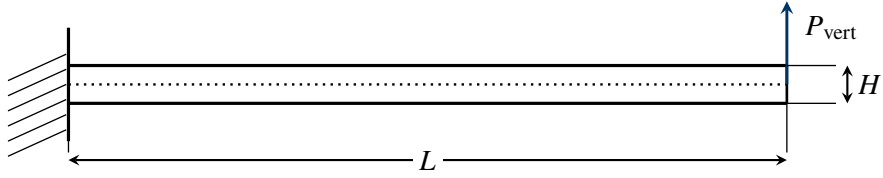
\begin{figure}[H]
\centering
\begin{tikzpicture}[x=1cm,y=1cm,>=stealth]
  \def\Lbeam{9.5}
  \def\Hbeam{0.5}

  \draw[line width=1.2pt] (0,\Hbeam/2) -- (\Lbeam,\Hbeam/2);
  \draw[line width=1.2pt] (0,-\Hbeam/2) -- (\Lbeam,-\Hbeam/2);
  \draw[dotted, line width=0.9pt] (0,0) -- (\Lbeam,0);
  \draw[line width=1.2pt] (0,-0.75) -- (0,0.75);
  \foreach \y in {-0.6,-0.4,...,0.6}{
    \draw[line width=0.5pt] (-0.03,\y) -- (-0.8,\y-0.35);
  }
  \draw[line width=1.0pt] (\Lbeam,-\Hbeam/2) -- (\Lbeam,\Hbeam/2);

  \draw[->,line width=1.1pt,SwanseaBlue]
    (\Lbeam,0) -- (\Lbeam,1.1);
  \node[right] at (\Lbeam+0.1,0.75) {$P_{\mathrm{vert}}$};

  \draw[<->,line width=0.8pt] (0,-1.0) -- (\Lbeam,-1.0);
  \node[fill=white,inner sep=2pt] at (\Lbeam/2,-1.0) {$L$};
  \draw[line width=0.5pt] (0,-\Hbeam/2) -- (0,-0.9);
  \draw[line width=0.5pt] (\Lbeam,-\Hbeam/2) -- (\Lbeam,-0.9);

  \draw[<->,line width=0.8pt]
    (\Lbeam+0.8,-\Hbeam/2) -- (\Lbeam+0.8,\Hbeam/2);
  \node[right] at (\Lbeam+0.8,0) {$H$};
  \draw[line width=0.5pt]
    (\Lbeam,\Hbeam/2) -- (\Lbeam+0.65,\Hbeam/2);
  \draw[line width=0.5pt]
    (\Lbeam,-\Hbeam/2) -- (\Lbeam+0.65,-\Hbeam/2);
\end{tikzpicture}
\caption{Cantilever beam subjected to the vertical tip load $P_{\mathrm{vert}}$. The beam has $L=100\,\mathrm{mm}$, $H=2\,\mathrm{mm}$, and out-of-plane width $W=1\,\mathrm{mm}$.}
\label{fig:tip_load_setup}
\end{figure}

The first loading case is the vertical force shown schematically in Figure~\ref{fig:tip_load_setup}. The point load at the free-end in the positive $Y$-direction is
\begin{align}
P_{\mathrm{vert}}
&=
\frac{4E_YI}{L^2}
=
1.616\times10^{-4}\,\mathrm{N}.
\label{eq:Pvert_def}
\end{align}

\begin{figure}[H]
    \centering
    \includegraphics[width=0.8\textwidth]{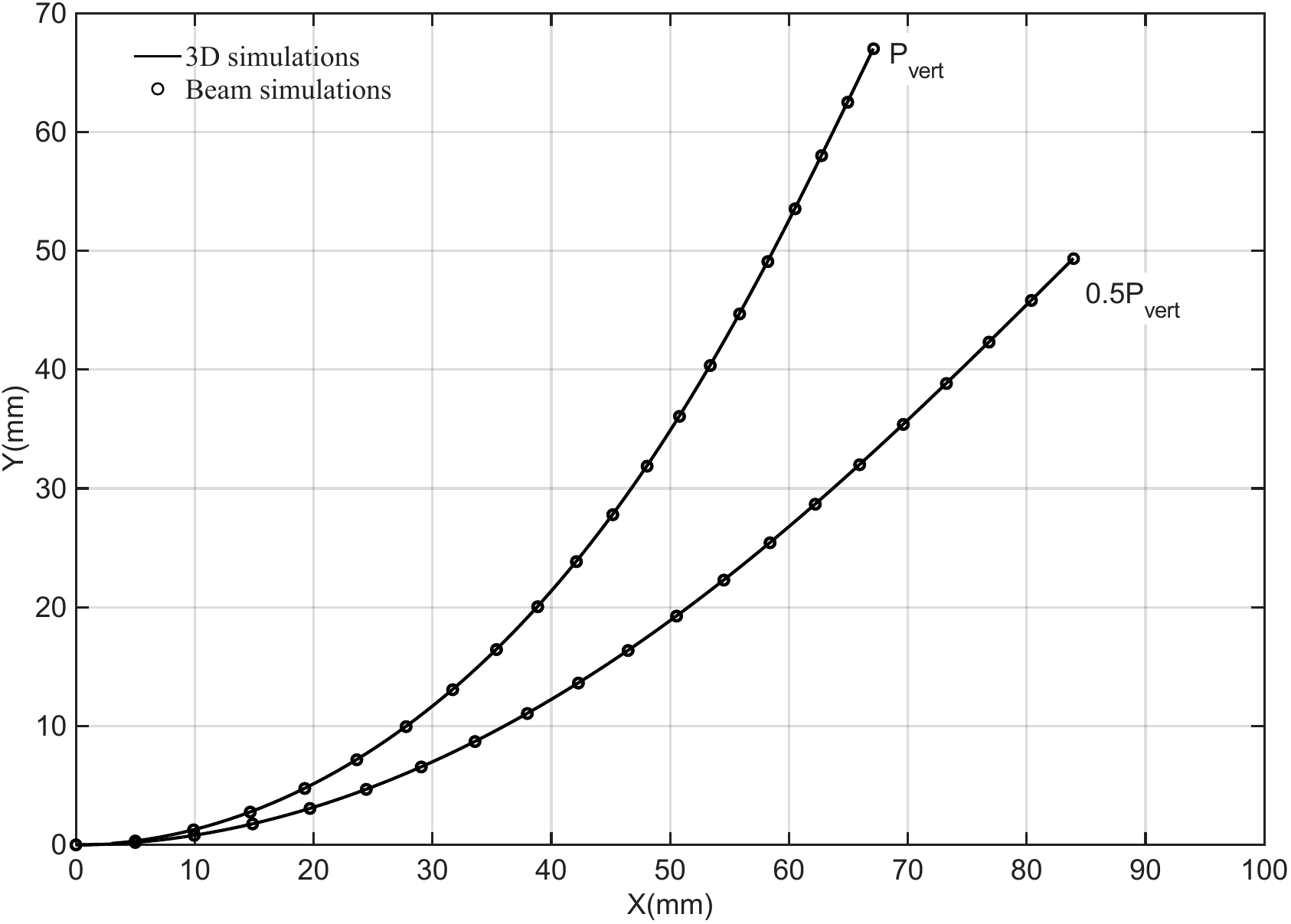}
    \caption{Deformed centrelines under the vertical tip load at 50\% and
100\% of $P_{\mathrm{vert}}$. Solid curves represent the centreline extracted from the three-dimensional continuum simulations and open circles represent the beam simulations. Coordinates are in millimetres.}
    \label{fig:linear_tip_load_comparison}
\end{figure}

The beam solution is compared with a three-dimensional continuum model discretised using quadratic hexahedral elements.
Figure~\ref{fig:linear_tip_load_comparison} compares the deformed centreline
of the beam model with that extracted from the three-dimensional continuum
simulation at $0.5P_{\mathrm{vert}}$ and $P_{\mathrm{vert}}$. This comparison shows that the geometrically exact beam reproduces the global three-dimensional deformation accurately.

\subsubsection{Cantilever beam with end moment}
\label{sec:linear_end_moment}

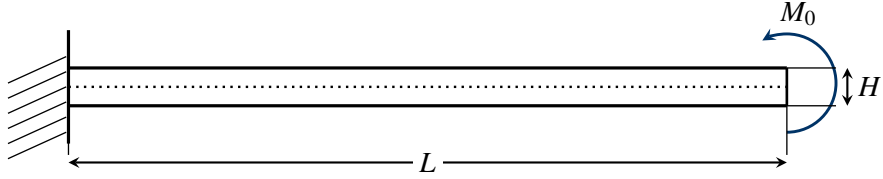
\begin{figure}[H]
\centering
\begin{tikzpicture}[x=1cm,y=1cm,>=stealth]
  \def\Lbeam{9.5}
  \def\Hbeam{0.5}

  \draw[line width=1.2pt] (0,\Hbeam/2) -- (\Lbeam,\Hbeam/2);
  \draw[line width=1.2pt] (0,-\Hbeam/2) -- (\Lbeam,-\Hbeam/2);
  \draw[dotted, line width=0.9pt] (0,0) -- (\Lbeam,0);
  \draw[line width=1.2pt] (0,-0.75) -- (0,0.75);
  \foreach \y in {-0.6,-0.4,...,0.6}{
    \draw[line width=0.5pt] (-0.03,\y) -- (-0.8,\y-0.35);
  }
  \draw[line width=1.0pt] (\Lbeam,-\Hbeam/2) -- (\Lbeam,\Hbeam/2);

  \draw[->,line width=1.1pt,SwanseaBlue]
    (\Lbeam+0.0,-0.6)
    arc[start angle=-90,end angle=120,radius=0.65];
  \node[above] at (\Lbeam+0.15,0.68) {$M_0$};

  \draw[<->,line width=0.8pt] (0,-1.0) -- (\Lbeam,-1.0);
  \node[fill=white,inner sep=2pt] at (\Lbeam/2,-1.0) {$L$};
  \draw[line width=0.5pt] (0,-\Hbeam/2) -- (0,-0.9);
  \draw[line width=0.5pt] (\Lbeam,-\Hbeam/2) -- (\Lbeam,-0.9);

  \draw[<->,line width=0.8pt]
    (\Lbeam+0.8,-\Hbeam/2) -- (\Lbeam+0.8,\Hbeam/2);
  \node[right] at (\Lbeam+0.8,0) {$H$};
  \draw[line width=0.5pt]
    (\Lbeam,\Hbeam/2) -- (\Lbeam+0.65,\Hbeam/2);
  \draw[line width=0.5pt]
    (\Lbeam,-\Hbeam/2) -- (\Lbeam+0.65,-\Hbeam/2);
\end{tikzpicture}
\caption{Cantilever beam subjected to the prescribed end moment $M_0$. The beam has $L=100\,\mathrm{mm}$, $H=2\,\mathrm{mm}$, and out-of-plane width $W=1\,\mathrm{mm}$.}
\label{fig:end_moment_setup}
\end{figure}

The second loading case is the prescribed end moment shown in Figure~\ref{fig:end_moment_setup}. The reference moment corresponding to the curvature of one complete circular configuration is
\begin{align}
M_0
&=
\frac{2\pi E_YI}{L}
=
2.5384068641\times10^{-2}\,\mathrm{N\,mm}.
\label{eq:M0_linear_benchmark}
\end{align}

For this example, the conventional section law gives the constant curvature $\kappa$ and the exact centreline as
\[
\kappa=M/(E_YI),
\qquad
x(X)=\frac{\sin(\kappa X)}{\kappa},
\qquad
y(X)=\frac{1-\cos(\kappa X)}{\kappa},
\qquad 0\le X\le L.
\]

\begin{figure}[H]
    \centering
    \includegraphics[width=0.8\textwidth]{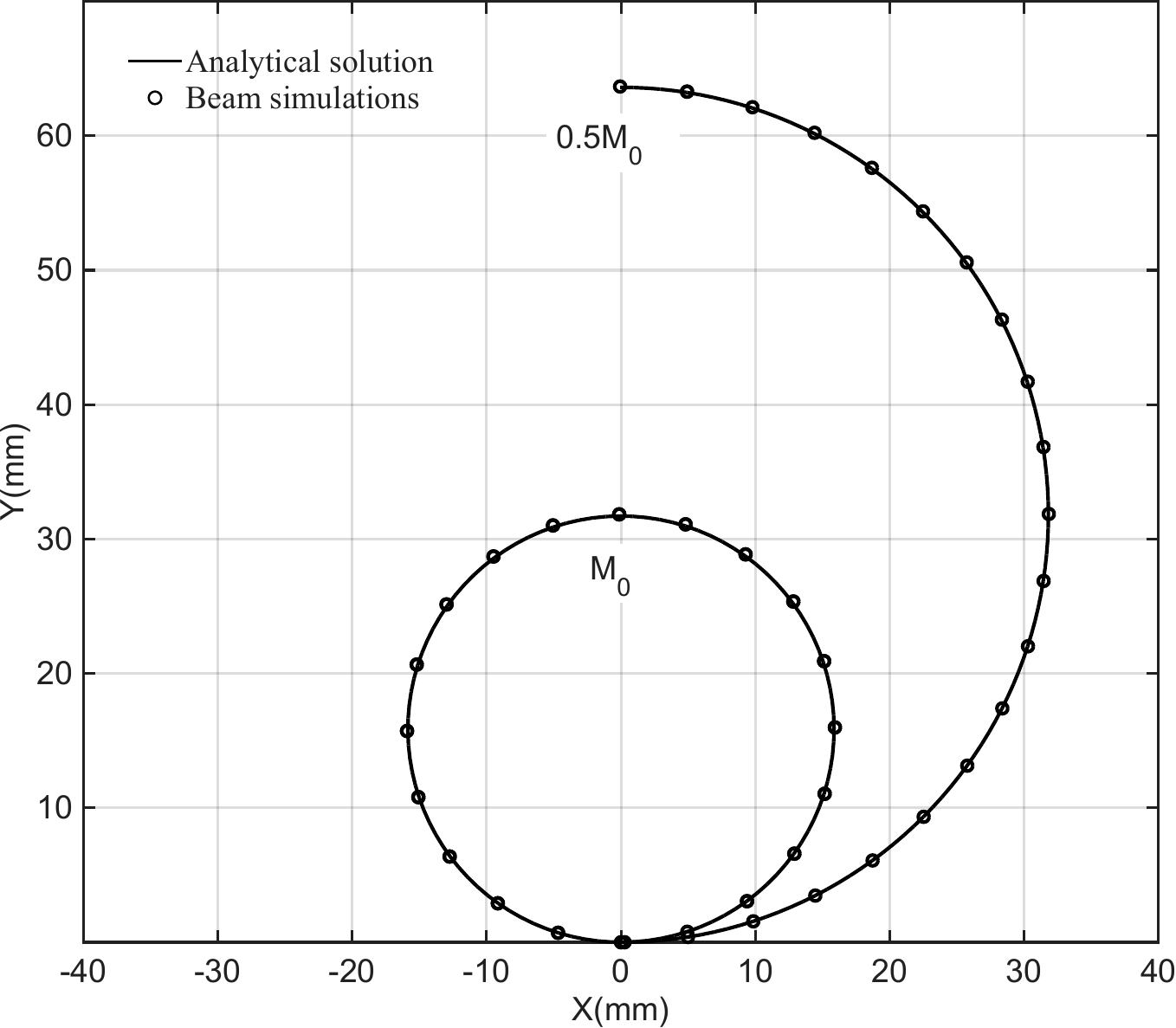}
    \caption{Deformed centrelines under the prescribed end moment at 50\% and 100\% of $M_0$. Solid curves represent the analytical pure-bending solutions and open circles represent the beam simulations. Coordinates are in millimetres.}
    \label{fig:linear_end_moment_comparison}
\end{figure}

The corresponding comparison is shown in Figure~\ref{fig:linear_end_moment_comparison}, for $0.5M_0$ and $M_0$. The beam points lie essentially on the analytical centrelines at both load levels. The close
agreement provides a direct verification of the finite-rotation kinematics and pure-bending response of the conventional beam formulation.

\FloatBarrier

\subsection*{Implications:}
\label{sec:linear_transition}


The comparisons in Figures~\ref{fig:linear_tip_load_comparison} and~\ref{fig:linear_end_moment_comparison} establish the accuracy of the conventional geometrically exact beam. 
The weak influence of Poisson's ratio observed in these comparisons is consistent with the conventional section law in Equation~\eqref{eq:H_linear}.
For a fixed Young's modulus and reference geometry, the axial and bending stiffnesses $E_YA$ and $E_YI$ remain unchanged, whereas Poisson's ratio
enters through the shear modulus $\mu$ in Equation~\eqref{eq:G_linear}, and therefore through the transverse shear stiffness $\kappa_s\mu A$.
Variations in $\nu$ consequently modify the shear response without directly affecting the axial or bending stiffnesses. Since both loading cases considered
for the present slender beam are predominantly bending dominated, the transverse shear contribution to the overall deformation is negligible. The resulting influence of $\nu$ on the global configuration is therefore limited,
consistent with the nearly coincident centrelines obtained in the numerical comparisons.

The above observations are specific to the linear material case considered here and should not be interpreted as a general independence from Poisson's ratio.
More importantly for the developments that follow, the conventional three-field kinematics do not admit an independent normal deformation of the cross-section.
The formulation with the thickness stretch included enriches the cross-sectional kinematics and treats the remaining out-of-plane constitutive response separately.

\FloatBarrier

\section{Enriched geometrically-exact beam formulation with thickness stretch}
\label{sec:hyperelastic_thickness}

The hyperelastic setting covers a broader range of constitutive responses than the conventional materially linear-elastic beam discussed in Section~\ref{sec:linear_beam}. For a three-dimensional hyperelastic solid, the deformation quantities retained in the beam model can still interact with transverse deformation through the material’s nonlinear response. As a result, if the adopted beam kinematics do not allow the cross-section to deform freely in the transverse direction, part of the constitutive behavior may be inadvertently restricted. To address this, the beam kinematics are enriched to allow deformation through the in-plane thickness, while the out-of-plane response is treated separately through the chosen constitutive reduction. Under plane strain, out-of-plane deformation is prevented, resulting in a reaction stress. Under plane stress, the out-of-plane stress is assumed to be zero and the associated deformation is determined locally from the three-dimensional constitutive law. This separated approach allows transverse deformation and Poisson-type constitutive effects to be represented consistently within a hyperelastic beam formulation.


For the enriched formulation, the mechanical field vector is extended to
\begin{equation}
\mathbf q_{\gamma}(X)
=
\begin{bmatrix}
u(X) & v(X) & \beta(X) & \gamma(X)
\end{bmatrix}^{\mathrm T}.
\label{eq:q_gamma}
\end{equation}
Here, $\gamma(X)$ is the global enrichment field governing deformation through the in-plane beam thickness in the reference $Y$-direction.

\subsection{Thickness-stretch beam kinematics}
\label{sec:thickness_kinematics}

The enrichment is expressed through the dimensionless through-thickness scaling measure
\begin{subequations}
\label{eq:thickness_mapping}
\begin{equation}
\lambda_t(X,Y)
=
1+Y\gamma(X),
\label{eq:lambda_t}
\end{equation}
so that $\gamma$ has dimensions of inverse length. The current signed transverse offset of a material point from the centreline is
\begin{equation}
\rho(X,Y)
=
Y\lambda_t(X,Y)
=
Y+Y^2\gamma(X).
\label{eq:rho}
\end{equation}
The associated derivatives are
\begin{equation}
\rho_{,X}
=
Y^2\gamma_{,X},
\qquad
\rho_{,Y}
=
1+2Y\gamma.
\label{eq:rho_derivatives}
\end{equation}
\end{subequations}
In particular, $\rho_{,Y}$ is the local transverse stretch associated with the $Y$-mapping, whereas $\lambda_t$ is the scaling measure used to report the surface response.

Accordingly, the in-plane placement is given by
\begin{subequations}
\label{eq:placement_stretch}
\begin{align}
x(X,Y)
&=
X+u(X)+\rho(X,Y)\sin\beta(X),
\label{eq:stretch_x}
\\
y(X,Y)
&=
v(X)+\rho(X,Y)\cos\beta(X).
\label{eq:stretch_y}
\end{align}
\end{subequations}
The in-plane deformation gradient is written as
\begin{subequations}
\label{eq:F_components_stretch}
\begin{equation}
\mathbf F_{2D}
=
\begin{bmatrix}
F_{11} & F_{12}\\
F_{21} & F_{22}
\end{bmatrix},
\label{eq:F2D_stretch}
\end{equation}
where
\begin{align}
F_{11}
&=
1+u_{,X}
+\rho_{,X}\sin\beta
+\rho\cos\beta\,\beta_{,X},
\label{eq:F11_stretch}
\\
F_{12}
&=
\rho_{,Y}\sin\beta,
\label{eq:F12_stretch}
\\
F_{21}
&=
v_{,X}
+\rho_{,X}\cos\beta
-\rho\sin\beta\,\beta_{,X},
\label{eq:F21_stretch}
\\
F_{22}
&=
\rho_{,Y}\cos\beta.
\label{eq:F22_stretch}
\end{align}
\end{subequations}
Substituting Equation~\eqref{eq:rho_derivatives} gives the equivalent explicit form
\begin{equation}
\begin{aligned}
F_{11}
&=
1+u_{,X}
+Y^2\gamma_{,X}\sin\beta
+\left(Y+Y^2\gamma\right)\cos\beta\,\beta_{,X},
\\
F_{12}
&=
(1+2Y\gamma)\sin\beta,
\\
F_{21}
&=
v_{,X}
+Y^2\gamma_{,X}\cos\beta
-\left(Y+Y^2\gamma\right)\sin\beta\,\beta_{,X},
\\
F_{22}
&=
(1+2Y\gamma)\cos\beta.
\end{aligned}
\end{equation}

The in-plane right Cauchy-Green tensor is
\begin{subequations}
\label{eq:strain_components_stretch}
\begin{equation}
\mathbf C_{2D}
=
\mathbf F_{2D}^{\mathrm T}\mathbf F_{2D},
\label{eq:C2D}
\end{equation}
and the Green-Lagrange strain components are
\begin{align}
E_{11}
&=
\frac12\left(F_{11}^2+F_{21}^2-1\right),
\label{eq:E11_stretch}
\\
E_{22}
&=
\frac12\left(F_{12}^2+F_{22}^2-1\right),
\label{eq:E22_stretch}
\\
E_{12}=E_{21}
&=
\frac12\left(F_{11}F_{12}+F_{21}F_{22}\right).
\label{eq:E12_stretch}
\end{align}
The two-dimensional Jacobian is
\begin{equation}
J_{2D}
=
\det\mathbf F_{2D}
=
F_{11}F_{22}-F_{12}F_{21}.
\label{eq:J2D_stretch}
\end{equation}
\end{subequations}

Therefore, defining the in-plane strain and stress vectors as
\begin{equation}
\mathbf e =
\begin{bmatrix}
E_{11}\\
E_{22}\\
2E_{12}
\end{bmatrix},
\qquad
\mathbf s =
\begin{bmatrix}
S_{11}\\
S_{22}\\
S_{12}
\end{bmatrix},
\label{eq:inplane_vectors}
\end{equation}
provides a common engineering-shear representation for both constitutive reductions. 
In the plane-strain and plane-stress formulations of
Sections~\ref{sec:plane_strain} and~\ref{sec:plane_stress}, respectively,
the superscripts $(\cdot)^{\mathrm{pe}}$ and $(\cdot)^{\mathrm{ps}}$
are used to distinguish quantities associated with the two constitutive
reductions. A barred quantity, $\bar{(\cdot)}$, denotes the corresponding
reduced form restricted to the retained in-plane description.

For the isotropic material specialisations used in the numerical examples and collected in \ref{app:material_specialisations}, the Lam\'e parameters are
\begin{equation}
\mu
=
\frac{E_Y}{2(1+\nu)},
\qquad
\lambda_{\mathrm L}
=
\frac{E_Y\nu}{(1+\nu)(1-2\nu)}.
\label{eq:Lame_common}
\end{equation}
Thus, for $\nu=0$, the first Lam\'e parameter satisfies
$\lambda_{\mathrm L}=0$, and the associated normal constitutive coupling
vanishes.

\subsection{Rationale for retaining the full strain and stress fields}
\label{sec:enriched_decomposition}

An important consequence of introducing the thickness-stretch field is that
the interpretation of the beam response in terms of separate axial, shear,
and bending modes becomes less direct. In the conventional three-field
formulation of Section~\ref{sec:linear_formulation}, the retained
through-thickness strain distribution is
$E_{11}=E_0+YK_b$, $2E_{12}=\Gamma$, and $E_{22}=0$. The quantities
$E_0$ and $\Gamma$ are independent of $Y$ and describe the centreline axial
and transverse shear deformation, respectively, whereas the term $YK_b$
provides the first-order bending contribution. Accordingly, the strain tensor
is represented by the three fixed generators introduced in
Equation~\eqref{eq:strain_generator_expansion},
\begin{equation*}
\mathbf G_{E_0}
=
\mathbf a_1\otimes\mathbf a_1,
\qquad
\mathbf G_{\Gamma}
=
\frac12
\left(
\mathbf a_1\otimes\mathbf a_2
+
\mathbf a_2\otimes\mathbf a_1
\right),
\qquad
\mathbf G_{K_b}
=
Y\,\mathbf a_1\otimes\mathbf a_1.
\end{equation*}
This simple structure leads directly to the axial-force, shear-force, and
bending-moment resultants in
Equations~\eqref{eq:Tp_general}--\eqref{eq:Mb_general}.
The introduction of the thickness-stretch field changes the
through-thickness structure of the strain. Using
$\rho=Y+Y^2\gamma$, $\rho_{,X}=Y^2\gamma_{,X}$, and
$\rho_{,Y}=1+2Y\gamma$ in the deformation gradient of
Equations~\eqref{eq:F11_stretch}--\eqref{eq:F22_stretch}, the exact in-plane
Green-Lagrange strain vector can be collected according to powers of $Y$ as
\begin{equation}
\mathbf e(X,Y)
=
\sum_{m=0}^{4}
Y^m\mathbf e^{(m)}(X),
\label{eq:enriched_Y_expansion}
\end{equation}
where
\begin{align*}
\mathbf e^{(0)}
&=
\begin{bmatrix}
E_0\\
0\\
\Gamma
\end{bmatrix},
&
\mathbf e^{(1)}
&=
\begin{bmatrix}
K_b\\
2\gamma\\
2\gamma\Gamma
\end{bmatrix},
\\[6pt]
\mathbf e^{(2)}
&=
\begin{bmatrix}
\gamma K_b
+\dfrac12\beta_{,X}^{\,2}
+\Gamma\gamma_{,X}\\
2\gamma^2\\
\gamma_{,X}
\end{bmatrix},
&
\mathbf e^{(3)}
&=
\begin{bmatrix}
\gamma\beta_{,X}^{\,2}\\
0\\
2\gamma\gamma_{,X}
\end{bmatrix},
\\[6pt]
\mathbf e^{(4)}
&=
\begin{bmatrix}
\dfrac12
\left(
\gamma^2\beta_{,X}^{\,2}
+
\gamma_{,X}^{\,2}
\right)\\
0\\
0
\end{bmatrix}.
\end{align*}
Here, the vectors $\mathbf e^{(m)}$ represent the coefficients linked to the corresponding powers of $Y$ and serve as weights, and the through-thickness polynomial basis is given by $1, Y, Y^2, Y^3,$ and $Y^4$.
In comparison, the conventional formulation
retains only the constant axial and shear modes together with the
first-order axial mode associated with bending.

To place the enriched representation in the same tensorial form as the
conventional beam, let $\mathbf a_1$ and $\mathbf a_2$ denote the Cartesian
basis vectors in the reference $X$- and $Y$-directions, as introduced in
Section~\ref{sec:linear_formulation}. For each polynomial order
$m=0,\ldots,4$, define the fixed through-thickness strain generators
\[
\mathbf G_{11}^{(m)}=Y^m\,\mathbf a_1\otimes\mathbf a_1,\qquad
\mathbf G_{22}^{(m)}=Y^m\,\mathbf a_2\otimes\mathbf a_2,\qquad
\mathbf G_{12}^{(m)}=\frac{Y^m}{2}\left(\mathbf a_1\otimes\mathbf a_2+\mathbf a_2\otimes\mathbf a_1\right).
\]
Writing $\mathbf e^{(m)} = \left[e_1^{(m)}\;\; e_2^{(m)}\;\; e_3^{(m)}\right]^{\mathrm T}$, 
where the third component is associated with the engineering-shear strain
$2E_{12}$, the retained in-plane Green-Lagrange strain tensor may be
expressed exactly as
\begin{equation}
\mathbf E(X,Y)
=
\sum_{m=0}^{4}
\left[
e_1^{(m)}\mathbf G_{11}^{(m)}
+
e_2^{(m)}\mathbf G_{22}^{(m)}
+
e_3^{(m)}\mathbf G_{12}^{(m)}
\right].
\label{eq:enriched_tensor_generator_expansion}
\end{equation}
The generators are fixed with respect
to the deformation. Their dependence on $Y$ only prescribes the corresponding
through-thickness mode. In this sense, the enriched formulation admits the
same type of generator construction as the conventional beam, but requires a
larger set of modes.

The relation with the conventional formulation is seen directly from
Equation~\eqref{eq:enriched_tensor_generator_expansion}. The conventional strain
tensor may be recovered from the subset
\begin{equation}
\mathbf E_{\mathrm{conv}}
=
E_0\mathbf G_{11}^{(0)}
+
\Gamma\mathbf G_{12}^{(0)}
+
K_b\mathbf G_{11}^{(1)}.
\label{eq:conventional_subset_enriched_basis}
\end{equation}
Thus, the conventional model retains a constant axial mode, a constant shear
mode, and a first-order bending mode. The thickness-stretch enrichment
introduces additional transverse-normal and shear contributions and extends
the axial, transverse-normal, and shear strain distributions to higher powers
of $Y$.

The stress quantities associated with these generators can be obtained in the
same variational manner as in the conventional formulation. For a symmetric
second Piola-Kirchhoff stress tensor $\mathbf S$, we define
\begin{align*}
R_{11}^{(m)}
&=
\int_A
\mathbf G_{11}^{(m)}:\mathbf S
\,\mathrm dA
=
\int_A
Y^mS_{11}\,\mathrm dA,
\\
R_{22}^{(m)}
&=
\int_A
\mathbf G_{22}^{(m)}:\mathbf S
\,\mathrm dA
=
\int_A
Y^mS_{22}\,\mathrm dA,
\\
R_{12}^{(m)}
&=
\int_A
\mathbf G_{12}^{(m)}:\mathbf S
\,\mathrm dA
=
\int_A
Y^mS_{12}\,\mathrm dA.
\end{align*}
For each order $m$, these quantities may be collected as $\mathbf{R}^{(m)} = \left[ R_{11}^{(m)}\;\; R_{22}^{(m)}\;\; R_{12}^{(m)} \right]^{\mathrm T}$.
The corresponding contribution to the internal virtual work is then
\begin{equation}
\delta\Pi_{\mathrm{int}}
=
\int_0^L
\sum_{m=0}^{4}
\delta\mathbf e^{(m)\mathrm T}
\mathbf R^{(m)}
\,\mathrm dX.
\label{eq:enriched_virtual_work_moments}
\end{equation}
Hence, fixed strain generators and their associated stress moments can be
defined for the enriched formulation in direct analogy with the conventional
beam.

The difference lies in the structure of the strain coefficients and in their
constitutive interpretation. For the conventional beam, $E_0$, $\Gamma$,
and $K_b$ can be identified directly with axial, shear, and bending
deformation. In the enriched formulation, this separation is lost even though
the tensor basis itself remains well defined. The first-order coefficient
$\mathbf e^{(1)}$ contains the bending term $K_b$, the transverse-normal
contribution $2\gamma$, and the coupled shear term $2\gamma\Gamma$.
At second order, the axial component contains
$\gamma K_b$, $\beta_{,X}^{\,2}/2$, and
$\Gamma\gamma_{,X}$, while the transverse-normal and shear components
contain $2\gamma^2$ and $\gamma_{,X}$, respectively. Higher-order terms
introduce further products of $\gamma$, $\gamma_{,X}$, and
$\beta_{,X}$. In particular, $\gamma_{,X}$ first appears at second order.
Restricting the enriched strain to zeroth- and first-order terms would
therefore remove part of the kinematics introduced by the thickness-stretch
field.

An additional distinction stems from the constitutive behaviour.
For a general hyperelastic material, the relation
$\mathbf S=\mathbf S(\mathbf E)$ is nonlinear. As a result, the stress
produced by Equation~\eqref{eq:enriched_tensor_generator_expansion}
cannot, in general, be decomposed into separate stress components
linked to the individual strain coefficients of
$\mathbf{e}^{(m)}$. Achieving such a decomposition would require a further
expansion and truncation of the constitutive relation and would be
specific to the chosen hyperelastic model. In the plane-stress setting,
the local computation of the out-of-plane stretch adds another
constitutive dependency at each material point.

For this reason, the generator representation above can be used to clarify the
through-thickness structure of the enriched strain field and its relation to
the conventional beam, rather than to introduce a separate set of section
constitutive variables. In the numerical formulation, the complete strain
tensor is retained at each cross-sectional quadrature point. The corresponding
stress tensor and consistent material tangent are evaluated directly from the
selected three-dimensional constitutive reduction, and their contributions to
the internal force and tangent are integrated numerically over the
cross-section. This avoids a further constitutive truncation while retaining
the complete kinematic coupling produced by the thickness-stretch
enrichment.

\subsection{Finite-element treatment of the thickness-stretch kinematics}
\label{sec:hyper_FE}

The enrichment adds one nodal field without changing the one-dimensional beam discretisation. The same shape functions are used for $u$, $v$, $\beta$, and $\gamma$, while the constitutive response is evaluated at cross-sectional quadrature points using the appropriate three-dimensional reduction.

Within an element,
\begin{subequations}
\label{eq:FE_kinematics_gamma}
\begin{equation}
\mathbf q_{\gamma}^{h}(X)
=
\sum_{a=1}^{n_e}N_a(X)\mathbf q_{\gamma a},
\qquad
\mathbf q_{\gamma a}
=
\begin{bmatrix}
u_a & v_a & \beta_a & \gamma_a
\end{bmatrix}^{\mathrm T},
\label{eq:FE_interpolation_gamma}
\end{equation}
with element vector
\begin{equation}
\mathbf q_{\gamma,e}
=
\begin{bmatrix}
u_1 & v_1 & \beta_1 & \gamma_1 & \cdots &
u_{n_e} & v_{n_e} & \beta_{n_e} & \gamma_{n_e}
\end{bmatrix}^{\mathrm T}.
\label{eq:q_gamma_element}
\end{equation}
The nonlinear strain-displacement matrix $\mathbf B_E$ is defined by
\begin{equation}
\delta\mathbf e
=
\mathbf B_E\,\delta\mathbf q_{\gamma,e}.
\label{eq:BE_definition}
\end{equation}
\end{subequations}
Its entries follow directly from the variations of
Equations~\eqref{eq:F11_stretch}-\eqref{eq:E12_stretch}. This definition is used
unchanged in the plane-strain and plane-stress material contributions below;
only the constitutive matrix supplied at the quadrature point differs.
The element residual and tangent contributions are evaluated by numerical
integration over the cross-section and along the beam axis. As in the
conventional beam formulation, reduced integration is employed in the axial
direction to alleviate shear locking in the slender-beam limit.

\subsection{Plane-strain formulation}
\label{sec:plane_strain}

The plane-strain model retains the enriched in-plane kinematics while suppressing deformation in the $Z$-direction. Accordingly, the out-of-plane strain vanishes and the corresponding second Piola-Kirchhoff stress is generally non-zero, acting as the reaction required to enforce the kinematic constraint.

\subsubsection{General hyperelastic material}

The constitutive definitions in Equation~\eqref{eq:constitutive_general_linear_section} are retained, now with the full three-dimensional strain tensor and indices $I,J,K,L\in\{1,2,3\}$. The reduction is written independently of a particular strain-energy function. The Saint-Venant-Kirchhoff specialisation used in the numerical examples is given in \ref{app:svk_plane_strain}.

Under plane strain,
\begin{subequations}
\label{eq:plane_strain_kinematics}
\begin{equation}
z(X,Y,Z)=Z,
\label{eq:z_plane_strain}
\end{equation}
and the three-dimensional deformation gradient is
\begin{equation}
\mathbf F^{\mathrm{pe}}
=
\begin{bmatrix}
F_{11} & F_{12} & 0\\
F_{21} & F_{22} & 0\\
0 & 0 & 1
\end{bmatrix}.
\label{eq:F_plane_strain}
\end{equation}
Hence,
\begin{equation}
F_{33}=1,
\qquad
J^{\mathrm{pe}}=J_{2D},
\label{eq:plane_strain_F33_J}
\end{equation}
and
\begin{equation}
\mathbf E^{\mathrm{pe}}
=
\begin{bmatrix}
E_{11} & E_{12} & 0\\
E_{12} & E_{22} & 0\\
0 & 0 & 0
\end{bmatrix},
\qquad
E_{33}=E_{13}=E_{23}=0.
\label{eq:E_plane_strain}
\end{equation}
\end{subequations}

The constitutive response is evaluated directly at
$E_{33}=E_{13}=E_{23}=0$. No restriction is imposed on the out-of-plane
stress, so
\begin{equation}
S_{33}\neq0
\qquad
\text{in general}.
\label{eq:S33_plane_strain}
\end{equation}

Because admissible increments satisfy
$D E_{33}=D E_{13}=D E_{23}=0$, the in-plane material tangent is the direct
in-plane block of the three-dimensional tangent,
\begin{subequations}
\label{eq:plane_strain_tangent_set}
\begin{equation}
L_{\alpha\beta\gamma\delta}^{\mathrm{pe}}
=
L_{\alpha\beta\gamma\delta}
\big|_{E_{33}=E_{13}=E_{23}=0},
\qquad
\alpha,\beta,\gamma,\delta\in\{1,2\}.
\label{eq:L_plane_strain}
\end{equation}
In the vector convention of Equation~\eqref{eq:inplane_vectors}, write
\begin{equation}
D\mathbf s
=
\mathbf C_{\mathrm{mat}}^{\mathrm{pe}}D\mathbf e,
\label{eq:Cmat_plane_strain_law}
\end{equation}
where
\begin{equation}
\mathbf C_{\mathrm{mat}}^{\mathrm{pe}}
=
\begin{bmatrix}
L_{1111} & L_{1122} & L_{1112}\\
L_{2211} & L_{2222} & L_{2212}\\
L_{1211} & L_{1222} & L_{1212}
\end{bmatrix}_{\mathrm{pe}}.
\label{eq:Cmat_plane_strain}
\end{equation}
\end{subequations}

The internal virtual work is
\begin{subequations}
\label{eq:plane_strain_element_contributions}
\begin{equation}
\delta\Pi_{\mathrm{int}}^{\mathrm{pe}}
=
\int_0^{L_e}\int_A
\mathbf s^{\mathrm T}\delta\mathbf e
\,\mathrm dA\,\mathrm dX.
\label{eq:virtual_work_plane_strain}
\end{equation}
Although $S_{33}$ is generally non-zero, it performs no virtual work because
$\delta E_{33}=0$. The corresponding material contribution to the element
tangent is
\begin{equation}
\mathbf K_{\mathrm{mat}}^{\mathrm{pe}}
=
\int_0^{L_e}\int_A
\mathbf B_E^{\mathrm T}
\mathbf C_{\mathrm{mat}}^{\mathrm{pe}}
\mathbf B_E
\,\mathrm dA\,\mathrm dX.
\label{eq:Kmat_plane_strain}
\end{equation}
\end{subequations}
The remaining geometric part of the element tangent follows from the
linearisation of the nonlinear strain-displacement relation, in the same
manner as for the conventional geometrically exact beam.

The Saint-Venant-Kirchhoff material used for the plane-strain benchmark is specialised from the above general reduction in \ref{app:svk_plane_strain}. That appendix gives the corresponding in-plane constitutive matrix, reduced strain-energy density, and out-of-plane reaction stress.

\subsection{Plane-stress formulation}
\label{sec:plane_stress}

The plane-stress formulation retains the same enriched in-plane kinematics but replaces the prescribed out-of-plane strain by the traction condition $S_{33}=0$. The corresponding stretch $\lambda_z$ is determined independently at each material quadrature point and is therefore an internal constitutive variable rather than a global beam field.

\subsubsection{General hyperelastic material}

The general reduction is written first so that the local solution for $\lambda_z$ and the consistent tangent condensation are independent of the chosen hyperelastic energy. At each material quadrature point, the in-plane deformation gradient is augmented locally as
\begin{subequations}
\label{eq:plane_stress_kinematics}
\begin{equation}
\mathbf F^{\mathrm{ps}}
=
\begin{bmatrix}
F_{11} & F_{12} & 0\\
F_{21} & F_{22} & 0\\
0 & 0 & \lambda_z
\end{bmatrix}.
\label{eq:F_plane_stress}
\end{equation}
Its determinant is
\begin{equation}
J
=
\det\mathbf F^{\mathrm{ps}}
=
\lambda_z J_{2D}.
\label{eq:J_plane_stress}
\end{equation}
The Green-Lagrange strain tensor is
\begin{equation}
\mathbf E^{\mathrm{ps}}
=
\begin{bmatrix}
E_{11} & E_{12} & 0\\
E_{12} & E_{22} & 0\\
0 & 0 & E_{33}
\end{bmatrix},
\qquad
E_{33}
=
\frac12(\lambda_z^2-1).
\label{eq:E_plane_stress}
\end{equation}
\end{subequations}

The plane-stress condition is
\begin{equation}
S_{33}=0.
\label{eq:plane_stress_condition}
\end{equation}
Note that this plane-stress constraint is distinct from exact incompressibility. The
former allows the out-of-plane stretch $\lambda_z$ to adjust so that
$S_{33}=0$, whereas exact incompressibility would additionally require $J=1$.
For a general hyperelastic material the plane-stress condition defines the local scalar equation
\begin{subequations}
\label{eq:plane_stress_local_solve}
\begin{equation}
g(\lambda_z)
=
S_{33}\!\left(E_{11},E_{12},E_{22},E_{33}(\lambda_z)\right)
=
0.
\label{eq:g_lambda}
\end{equation}
Since
\begin{equation}
\frac{\partial E_{33}}{\partial\lambda_z}
=
\lambda_z,
\end{equation}
the derivative of the local residual is
\begin{equation}
\frac{\partial S_{33}}{\partial\lambda_z}
=
L_{3333}\lambda_z.
\label{eq:dS33_dlambda}
\end{equation}
A local Newton update can therefore be written as
\begin{equation}
\lambda_z^{(k+1)}
=
\lambda_z^{(k)}
-
\frac{S_{33}^{(k)}}{\lambda_z^{(k)}L_{3333}^{(k)}}.
\label{eq:lambda_Newton_general}
\end{equation}
\end{subequations}

The in-plane material tangent must be reduced consistently with the same
constraint. Linearisation of $S_{33}=0$ gives
\begin{subequations}
\label{eq:plane_stress_condensation}
\begin{equation}
D E_{33}
=
-
\frac{L_{33\alpha\beta}}{L_{3333}}
D E_{\alpha\beta},
\label{eq:DE33_condensed}
\end{equation}
and consequently
\begin{equation}
D S_{\alpha\beta}
=
\bar L_{\alpha\beta\gamma\delta}^{\mathrm{ps}}
D E_{\gamma\delta},
\label{eq:DS_condensed}
\end{equation}
where
\begin{equation}
\bar L_{\alpha\beta\gamma\delta}^{\mathrm{ps}}
=
L_{\alpha\beta\gamma\delta}
-
\frac{L_{\alpha\beta33}L_{33\gamma\delta}}{L_{3333}},
\qquad
\alpha,\beta,\gamma,\delta\in\{1,2\}.
\label{eq:L_condensed}
\end{equation}
\end{subequations}

Equation~\eqref{eq:L_condensed} has the structure of a Schur complement of the
three-dimensional material tangent with respect to the out-of-plane component.
Indeed, the plane-stress constraint couples the retained in-plane strain
increments to the eliminated increment $D E_{33}$ through the linearised
condition $D S_{33}=0$. Eliminating $D E_{33}$ using
Equation~\eqref{eq:DE33_condensed} therefore modifies the direct in-plane tangent by
the coupling through the out-of-plane direction, giving precisely the
condensed tensor in Equation~\eqref{eq:L_condensed}. The resulting tangent is thus
not an independently prescribed two-dimensional constitutive law, but the
consistent in-plane tangent of the original three-dimensional hyperelastic
material restricted to the plane-stress constraint. In this way, the local
adjustment of the out-of-plane deformation is transferred consistently to the
global beam equations without introducing $\lambda_z$ as an additional global
unknown.

In the vector notation of Equation~\eqref{eq:inplane_vectors},
\begin{subequations}
\label{eq:plane_stress_matrix_set}
\begin{equation}
D\mathbf s
=
\bar{\mathbf C}_{\mathrm{mat}}^{\mathrm{ps}}D\mathbf e,
\label{eq:C_plane_stress_law}
\end{equation}
with
\begin{equation}
\bar{\mathbf C}_{\mathrm{mat}}^{\mathrm{ps}}
=
\renewcommand{\arraystretch}{1.3}
\begin{bmatrix}
\bar L_{1111}^{\mathrm{ps}} & \bar L_{1122}^{\mathrm{ps}} & \bar L_{1112}^{\mathrm{ps}} \\
\bar L_{2211}^{\mathrm{ps}} & \bar L_{2222}^{\mathrm{ps}} & \bar L_{2212}^{\mathrm{ps}} \\
\bar L_{1211}^{\mathrm{ps}} & \bar L_{1222}^{\mathrm{ps}} & \bar L_{1212}^{\mathrm{ps}}
\end{bmatrix}.
\label{eq:C_plane_stress}
\end{equation}
\end{subequations}

Because $S_{33}=0$, the internal virtual work contains only the in-plane
components,
\begin{subequations}
\label{eq:plane_stress_element_contributions}
\begin{equation}
\delta\Pi_{\mathrm{int}}^{\mathrm{ps}}
=
\int_0^{L_e}\int_A
\mathbf s^{\mathrm T}\delta\mathbf e
\,\mathrm dA\,\mathrm dX,
\label{eq:virtual_work_plane_stress}
\end{equation}
and the material contribution to the element tangent is
\begin{equation}
\mathbf K_{\mathrm{mat}}^{\mathrm{ps}}
=
\int_0^{L_e}\int_A
\mathbf B_E^{\mathrm T}
\bar{\mathbf C}_{\mathrm{mat}}^{\mathrm{ps}}
\mathbf B_E
\,\mathrm dA\,\mathrm dX.
\label{eq:Kmat_plane_stress}
\end{equation}
\end{subequations}
Hence, the distinction from plane strain is confined to the third
direction: plane strain prescribes $E_{33}=0$, whereas plane stress determines
$E_{33}$ locally through $S_{33}=0$.

The material-specific plane-stress reductions used in the numerical examples are collected in \ref{app:material_specialisations}. The Saint-Venant-Kirchhoff model, for which the plane-stress constraint can be eliminated analytically, is given in \ref{app:svk_plane_stress}. The compressible Neo-Hookean model, for which the out-of-plane stretch is obtained from a local nonlinear scalar equation, is given in \ref{app:nh_plane_stress}.

\subsection{Numerical examples}
\label{sec:hyperelastic_results}


The hyperelastic calculations consider the same cantilever end-moment benchmark introduced in Section~\ref{sec:linear_end_moment}. The shared computational settings and convergence behaviour are as stated in Section~\ref{sec:linear_numerical_examples}.

For both the linear-elastic and hyperelastic beam models considered in this
work, the primary quantity of interest is the global deformation of the beam
centreline. In the hyperelastic formulation, however, the additional local
cross-sectional quantities provide information on how the transverse
deformation adjusts to the underlying constitutive coupling, particularly as
the material approaches near-incompressible behaviour. This aspect is
relevant for many rubber-like and soft polymeric materials encountered in
experiments and practical applications, which are commonly represented as
nearly incompressible hyperelastic solids
\cite{Ogden1972Incompressible,SteinmannHossainPossart2012}. The local
quantities are therefore examined primarily to clarify the transverse
deformations accompanying the global beam response, rather than as independent
measures of the overall structural deformation.

Accordingly, the global response is assessed via the deformed centreline, whereas the local in-plane thickness behaviour is analysed using the surface measures $\lambda_t(X,+H/2)$ and $\lambda_t(X,-H/2)$ introduced in Equation~\eqref{eq:lambda_t}. For plane-stress conditions, the associated out-of-plane stretch $\lambda_z$ is also computed at the material integration points to quantify the local transverse relaxation imposed by the plane-stress constraint. Four quadrature points are employed across the beam thickness after conducting studies with different numbers of quadrature points; five or more points yield indistinguishable results for both the overall deformation and the reported local measures such as $\lambda_t$. 

\subsubsection{Plane-strain response and limitations}
\label{sec:results_plane_strain}

\begin{figure}[htbp]
    \centering

    \begin{subfigure}[t]{0.75\textwidth}
        \centering
        \includegraphics[width=\linewidth]
        {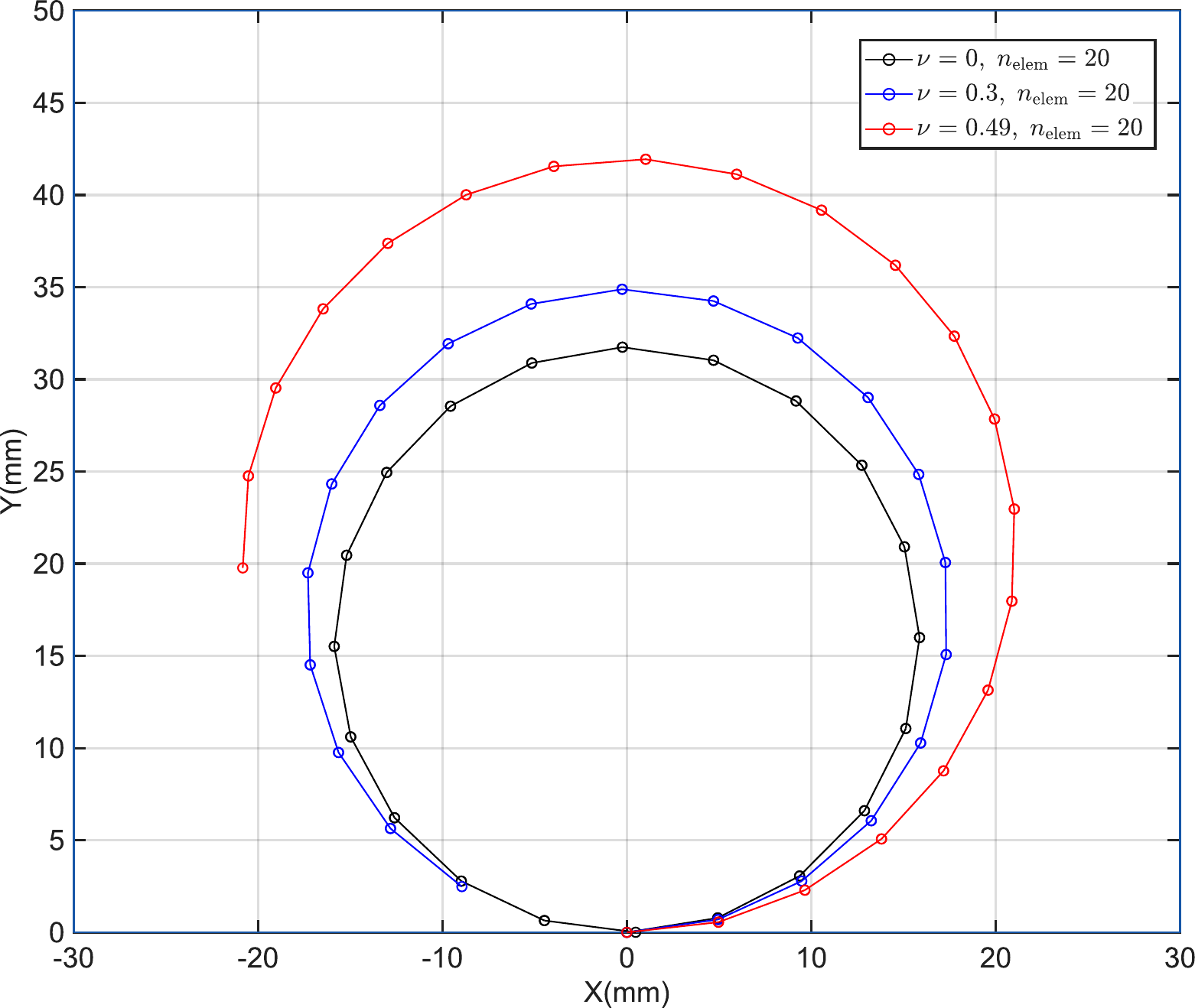}
        \caption{Influence of Poisson's ratio for
        $n_{\mathrm{elem}}=20$.}
        \label{fig:svk_plane_strain_def}
    \end{subfigure}

    \vspace{0.8em}

    \begin{subfigure}[t]{0.75\textwidth}
        \centering
        \includegraphics[width=\linewidth]
        {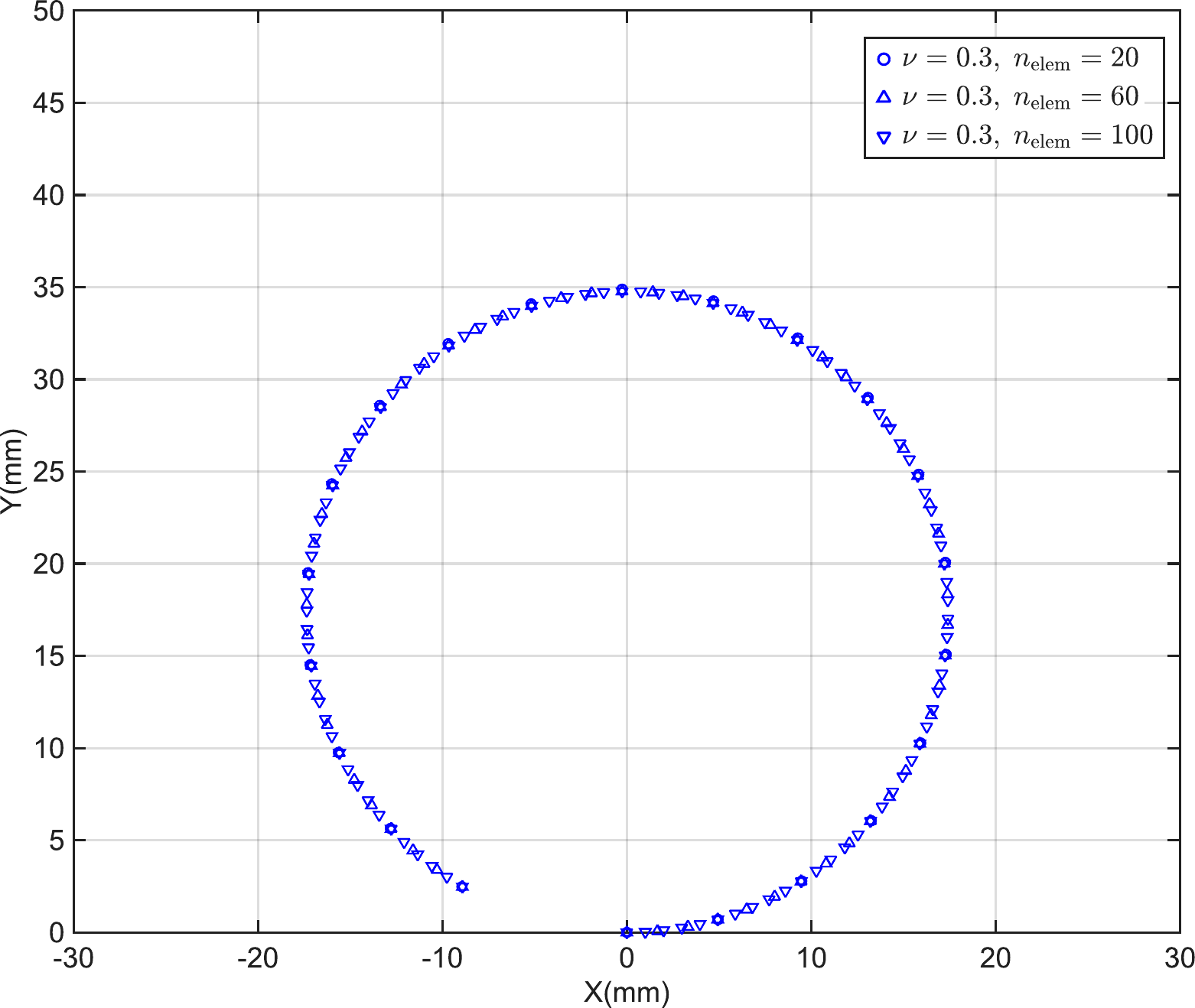}
        \caption{Mesh refinement for $\nu=0.3$.}
        \label{fig:svk_plane_strain_mesh}
    \end{subfigure}

    \caption{Deformed centreline of the Saint-Venant-Kirchhoff beam under
    plane strain for the cantilever end-moment problem of
    Section~\ref{sec:linear_end_moment}:
    (a) response for $\nu=0$, $0.3$, and $0.49$ using
    $n_{\mathrm{elem}}=20$; and
    (b) mesh-refinement response for $\nu=0.3$ using
    $n_{\mathrm{elem}}=20$, $60$, and $100$.}
    \label{fig:svk_plane_strain_global}
\end{figure}

\paragraph{Reduced Saint-Venant-Kirchhoff model}
The Saint-Venant-Kirchhoff plane-strain response, using the material specialisation in \ref{app:svk_plane_strain}, is evaluated for $\nu=0$, $0.3$, and $0.49$. 
Figure~\ref{fig:svk_plane_strain_global} shows the global centreline response.
Figure~\ref{fig:svk_plane_strain_def} compares the three Poisson-ratio cases using $n_{\mathrm{elem}}=20$, while Figure~\ref{fig:svk_plane_strain_mesh} provides a mesh-convergence check for $\nu=0.3$ using $n_{\mathrm{elem}}=20$, $60$, and $100$.

The $\nu=0$ solution in Figure~\ref{fig:svk_plane_strain_def} recovers the circular configuration associated with the corresponding end-moment benchmark, demonstrating the consistency with the conventional beam element.
As Poisson's ratio increases, the curvature obtained under the same prescribed moment decreases and the centreline becomes progressively more open.
The change is already apparent for $\nu=0.3$ and becomes considerably stronger for $\nu=0.49$.
This trend is consistent with the Poisson-dependent coupling in Equation~\eqref{eq:SVK_plane_strain_matrix} and the associated out-of-plane reaction in Equation~\eqref{eq:S33_SVK_plane_strain}.
The centreline comparison in Figure~\ref{fig:svk_plane_strain_mesh} shows satisfactory mesh convergence of the global response. The configurations obtained with $20$, $60$, and $100$ elements are practically coincident, with
no appreciable change in the overall curvature or free-end position. The global deformation is therefore already well represented by the discretisations considered here. More importantly, the mesh convergence reveals that the deformation does not change with mesh refinement, indicating that the difference from the unconstrained beam response concerns the modelling constraint rather than a simple sensitivity to the mesh employed.

\begin{figure}[p]
    \centering

    \begin{subfigure}[t]{0.75\textwidth}
        \centering
        \includegraphics[width=\linewidth]
        {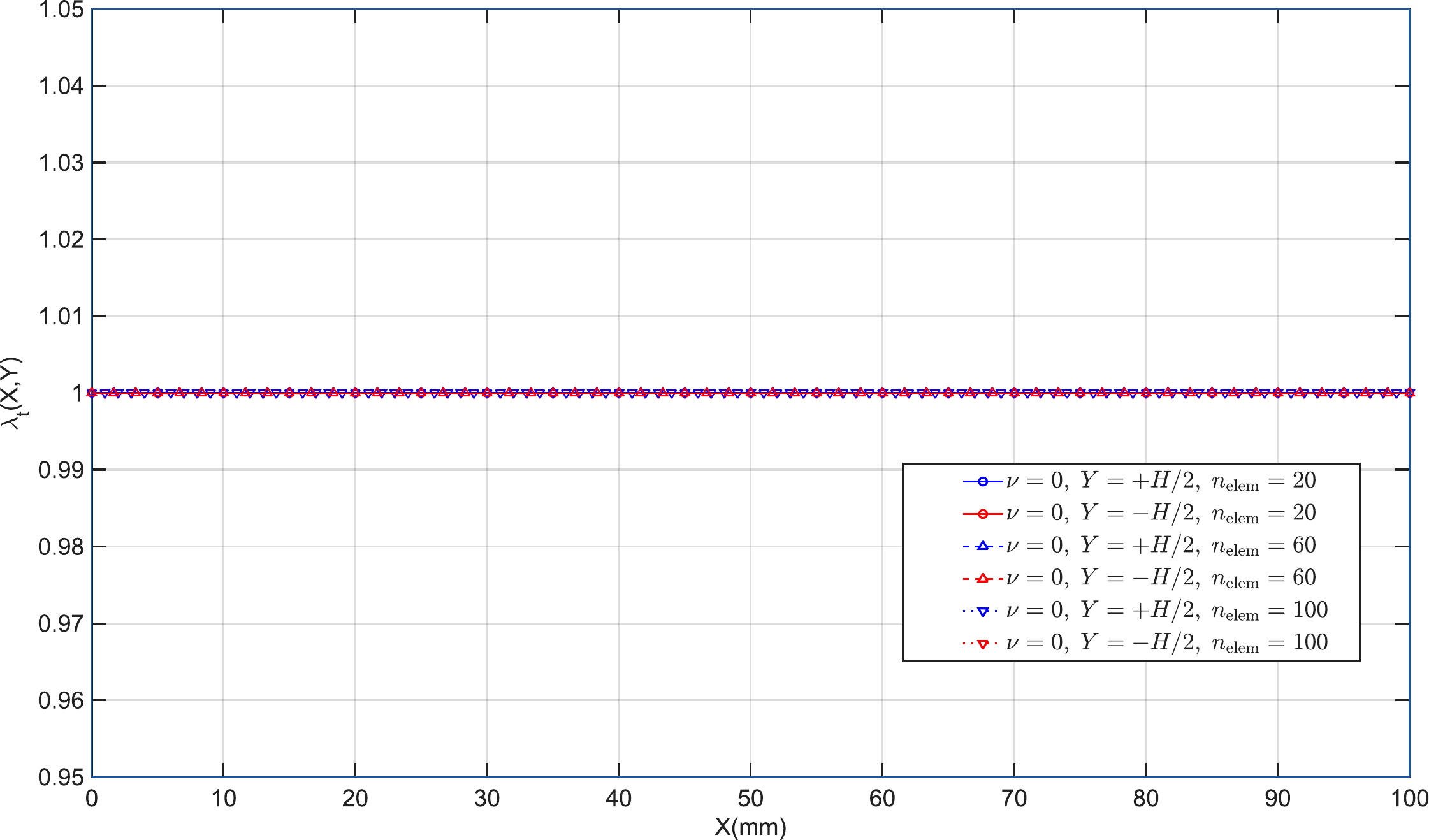}
        \caption{$\nu=0$.}
        \label{fig:svk_plane_strain_lambdat_nu0}
    \end{subfigure}

    \vspace{0.8em}

    \begin{subfigure}[t]{0.75\textwidth}
        \centering
        \includegraphics[width=\linewidth]
        {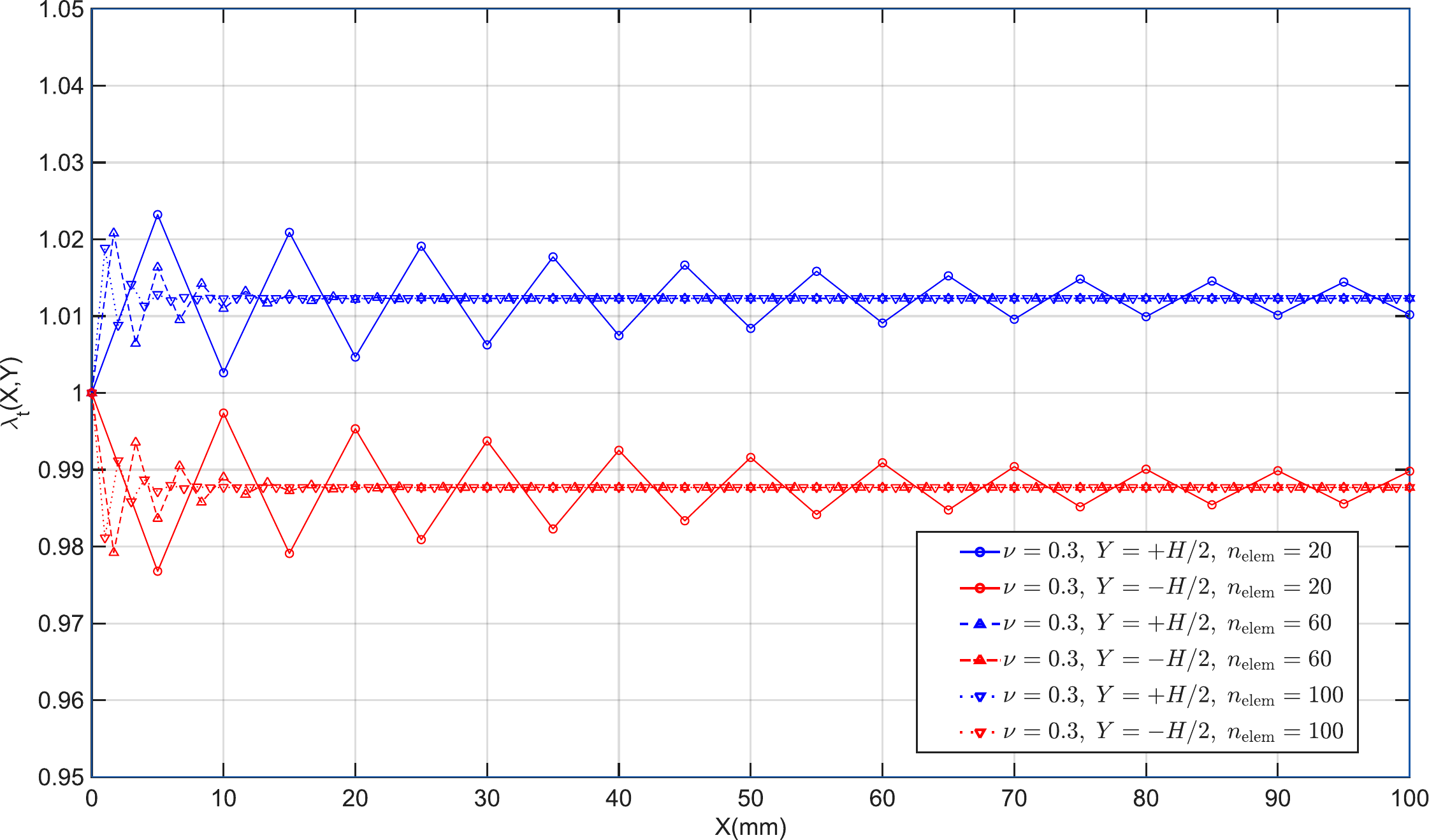}
        \caption{$\nu=0.3$.}
        \label{fig:svk_plane_strain_lambdat_nu03}
    \end{subfigure}

    \vspace{0.8em}

    \begin{subfigure}[t]{0.75\textwidth}
        \centering
        \includegraphics[width=\linewidth]
        {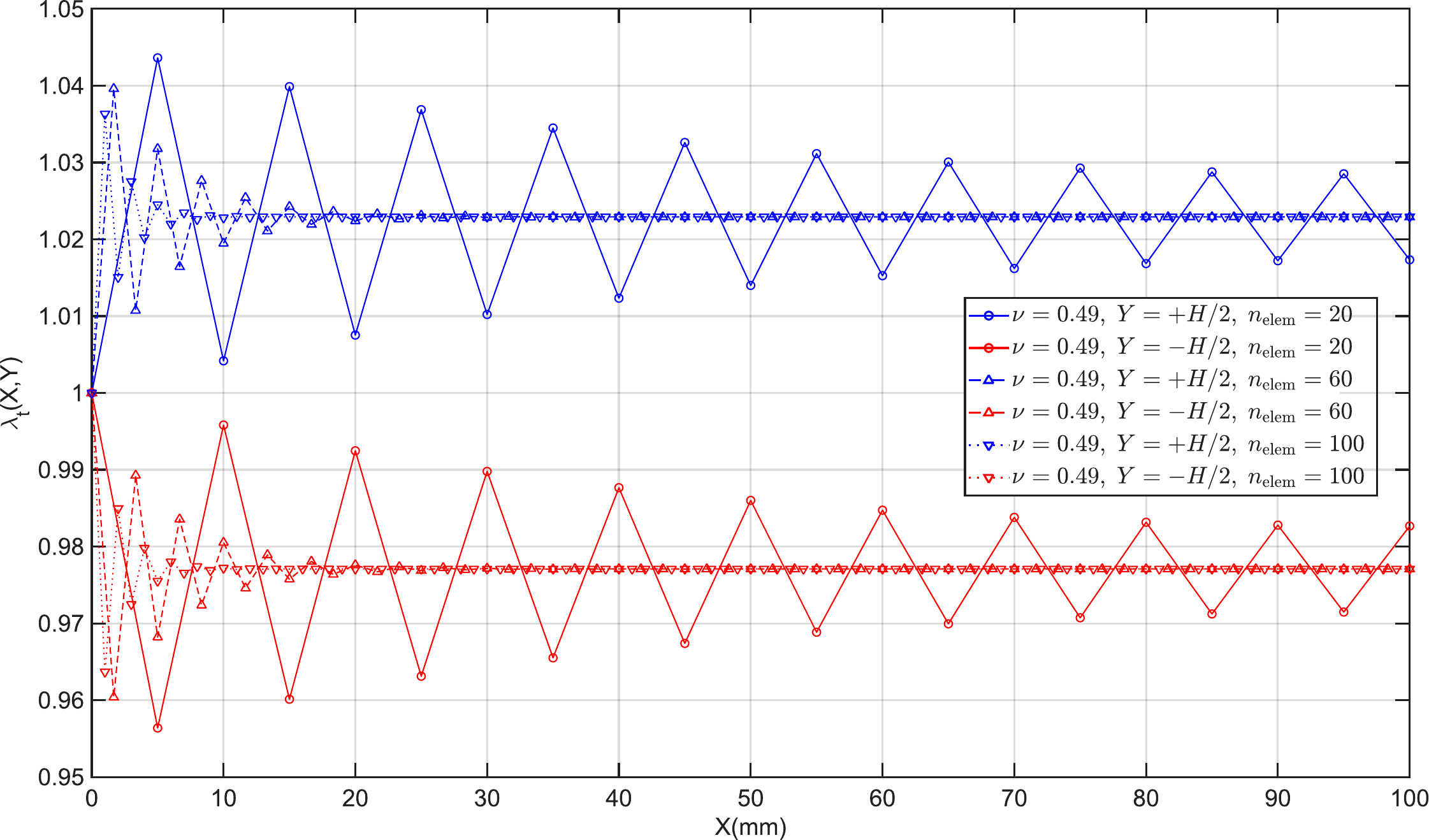}
        \caption{$\nu=0.49$.}
        \label{fig:svk_plane_strain_lambdat_nu049}
    \end{subfigure}

    \caption{Surface values of the through-thickness scaling measure
    $\lambda_t$ for the Saint-Venant-Kirchhoff beam under plane strain:
    (a) $\nu=0$, (b) $\nu=0.3$, and (c) $\nu=0.49$.
    Blue and red curves correspond to $Y=+H/2$ and $Y=-H/2$,
    respectively. Open circles, upward triangles, and downward triangles
    denote $n_{\mathrm{elem}}=20$, $60$, and $100$, respectively.}
    \label{fig:svk_plane_strain_lambdat}
\end{figure}

The enriched cross-sectional response is examined separately in Figure~\ref{fig:svk_plane_strain_lambdat} through the surface values of $\lambda_t$.
Results are shown using $n_{\mathrm{elem}}=20$, $60$, and $100$ to resolve the spatial structure of the response near the clamped end.
Blue and red curves correspond to $Y=+H/2$ and $Y=-H/2$, respectively, while the marker type identifies the number of elements.
For $\nu=0$, the surface values in Figure~\ref{fig:svk_plane_strain_lambdat_nu0} remain at $\lambda_t=1$ throughout the beam for all three discretisations, consistent with the absence of normal constitutive coupling in this limiting case.
For higher Poisson's ratio values, Figures~\ref{fig:svk_plane_strain_lambdat_nu03} and~\ref{fig:svk_plane_strain_lambdat_nu049} show a separation of the upper- and lower-surface responses, with the departure from unity increasing as the Poisson coupling becomes stronger.
In the finer-mesh solutions, the surface values approach nearly uniform interior levels away from the clamp. The axial adjustment near the constrained boundary should therefore be distinguished from the non-unit surface values in the beam interior. Smaller variations towards the loaded end are distinct from this constrained-boundary response.

Localised boundary-layer effects near constrained beam ends can influence local stress and cross-sectional quantities much more strongly than the overall deflection \cite{WangShi2012BoundaryLayer,Horgan1996SaintVenant}.
The importance of the boundary treatment for beams with deformable cross-sections has likewise been emphasised in geometrically and materially nonlinear formulations with extensible directors \cite{ChoiSauerKlinkel2021}.
Related nonlinear shell formulations retain through-the-thickness stretching as an additional kinematic contribution when using three-dimensional constitutive laws at large deformations \cite{BrankKorelcIbrahimbegovic2002}.
Shell boundary-layer studies provide further context for the distinction between local and interior responses. In thin conical shells, local compatibility bending near boundaries and changes in wall thickness supplements the interior membrane response \cite{BoyezSadowskiIzzuddin2018BoundaryLayer}. Analyses of shallow-shell models similarly identify boundary and interior layers whose characteristic length scales depend on the shell geometry and thickness \cite{PitkarantaMatacheSchwab2001Layers}. These studies illustrate how satisfying boundary constraints can require local deformation fields with spatial scales distinct from those of the interior response.

The different axial discretisations show how the local cross-sectional field is resolved. With $n_{\mathrm{elem}}=20$, the oscillatory variation extends over a larger portion of the discrete distribution. Increasing the number of elements produces flatter interior branches and a more localised adjustment near the clamped boundary, while the interior surface levels remain essentially unchanged. This local mesh sensitivity should be distinguished from the nearly coincident global centrelines in Figure~\ref{fig:svk_plane_strain_mesh}. The distinction is especially apparent for $\nu=0.49$, where the stronger Poisson coupling produces a larger local thickness response.

The need to resolve local cross-sectional effects has also been reported for higher-order beam finite elements, where local refinement captures clamped-end and load-discontinuity effects associated with cross-sectional warping \cite{GengDuanLi2017HigherOrder}. In plate and shell computations, Miazio and Zboi\'nski \cite{MiazioZboinski2023BoundaryLayer} detect boundary layers and resolve them using graded mesh refinement towards the boundary. Hakula \cite{Hakula2022HarmonicBoundaryLayers} likewise examines the resolution of boundary-layer scales using harmonic-extension finite elements, including a thin cylindrical-shell benchmark. These results support examining the resolution of boundary-sensitive fields separately from the global deformation.

This refinement does not remove the limitations of the beam kinematics: three-dimensional effects may become important near the clamp within a region whose characteristic extent scales with the beam thickness, so a beam model may not fully represent the local stress and deformation.
The same centreline and $\lambda_t$ results are obtained using either four or five Gauss-Legendre integration points through the thickness, confirming insensitivity to the two through-thickness quadrature rules for the axial discretisations considered here. This quadrature comparison concerns integration through the thickness, whereas the spatial resolution of the local field is assessed through the axial mesh comparison.

\subsubsection{Plane-stress response}
\label{sec:results_plane_stress}

\begin{figure}[htbp]
    \centering

    \begin{subfigure}[t]{0.8\textwidth}
        \centering
        \includegraphics[width=\linewidth]
        {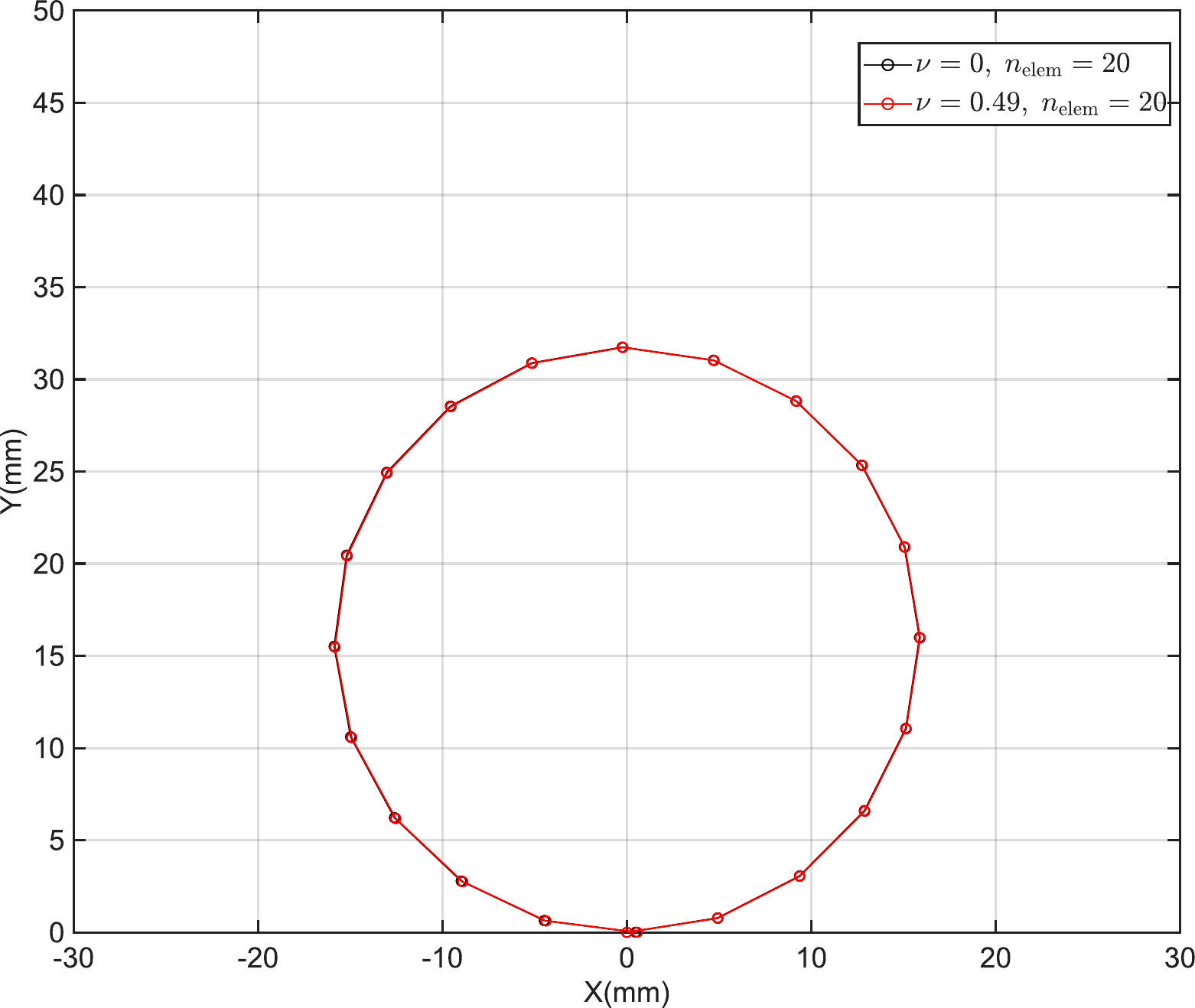}
        \caption{Deformed centreline for $\nu=0$ and $\nu=0.49$ using
        $n_{\mathrm{elem}}=20$.}
        \label{fig:svk_plane_stress_def}
    \end{subfigure}

    \vspace{0.8em}

    \begin{subfigure}[t]{0.85\textwidth}
        \centering
        \includegraphics[width=\linewidth]
        {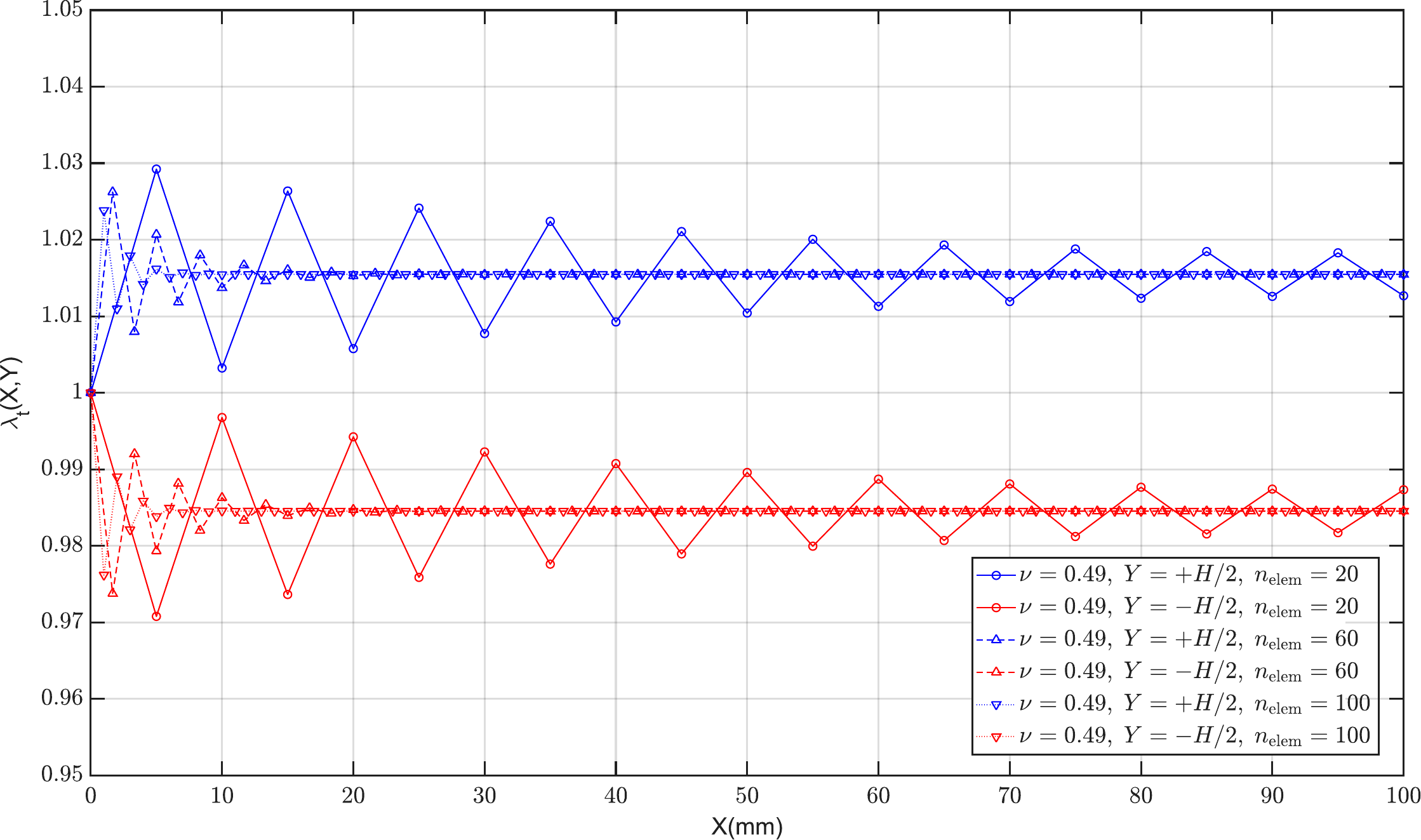}
        \caption{Surface values $\lambda_t(X,\pm H/2)$ for $\nu=0.49$
        using $n_{\mathrm{elem}}=20$, $60$, and $100$.}
        \label{fig:svk_plane_stress_lambda_t}
    \end{subfigure}

    \caption{Reduced Saint-Venant-Kirchhoff response under plane stress for
    the cantilever end-moment benchmark:
    (a) deformed centreline for $\nu=0$ and $\nu=0.49$, and
    (b) surface values of the through-thickness scaling measure $\lambda_t$
    for $\nu=0.49$. In (b), blue and red curves correspond to $Y=+H/2$ and
    $Y=-H/2$, respectively. Open circles, upward triangles, and downward
    triangles denote $n_{\mathrm{elem}}=20$, $60$, and $100$, respectively.}
    \label{fig:svk_plane_stress_results}
\end{figure}

\paragraph{Reduced Saint-Venant-Kirchhoff model}
The reduced Saint-Venant-Kirchhoff plane-stress calculations use the
specialisation in \ref{app:svk_plane_stress}. The same cantilever end-moment
benchmark is evaluated for $\nu=0$ and $\nu=0.49$.
Figure~\ref{fig:svk_plane_stress_results} compares the global deformation and
the in-plane thickness response.

The two centreline configurations in Figure~\ref{fig:svk_plane_stress_def} coincide over the complete beam.
Thus, in contrast to the plane-strain response in Figure~\ref{fig:svk_plane_strain_def}, changing Poisson's ratio from $\nu=0$ to $\nu=0.49$ produces no visible change in the global configuration for this plane-stress benchmark.
The representative $20$-element result is shown because the preceding plane-strain study has already established satisfactory centreline mesh convergence, and the same behaviour is obtained here.

The local thickness response remains non-trivial for $\nu=0.49$. Figure~\ref{fig:svk_plane_stress_lambda_t} shows that the upper- and lower-surface values separate about unity and approach nearly uniform levels
away from the constrained end of the beam.
The $20$-, $60$-, and $100$-element curves display the same behaviour near the constrained boundary discussed for plane strain: refinement resolves the short boundary adjustment more locally while leaving the interior response essentially unchanged. The finer distributions therefore characterise the local cross-sectional field rather than a change in the already established global bending configuration.

\begin{figure}[htbp]
    \centering
    \includegraphics[width=1\textwidth]
    {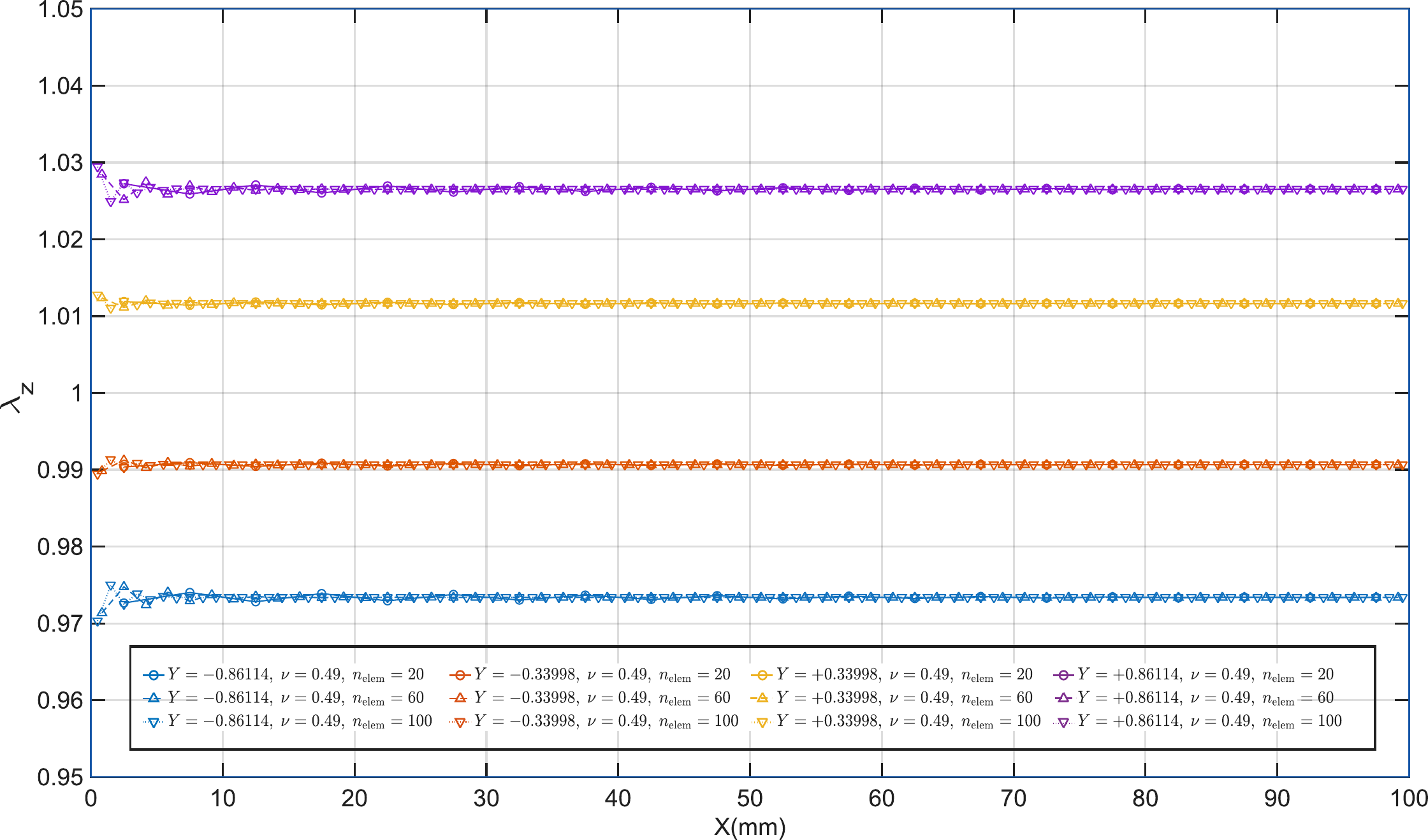}
    \caption{Out-of-plane stretch $\lambda_z$ for the reduced
    Saint-Venant-Kirchhoff plane-stress model with $\nu=0.49$. Values are
    evaluated at the axial integration point of each element and at the four
    through-thickness Gauss-Legendre locations. Blue, orange, yellow, and
    purple curves correspond, respectively, to
    $Y=-0.86114$, $-0.33998$, $+0.33998$, and $+0.86114\,\mathrm{mm}$.
    Open circles, upward triangles, and downward triangles denote
    $n_{\mathrm{elem}}=20$, $60$, and $100$, respectively.}
    \label{fig:svk_plane_stress_lambda_z}
\end{figure}

The out-of-plane response is examined through $\lambda_z$. Unlike $\lambda_t$, which follows from the global enrichment field, $\lambda_z$ is determined locally from the condition $S_{33}=0$ at the material quadrature
points. With one axial integration point in each uniform beam element, and four through-thickness quadrature points at every axial integration point, $20$-, $60$-, and $100$-element meshes thus
provide $80$, $240$, and $400$ local values, respectively, as shown in Figure \ref{fig:svk_plane_stress_lambda_z}.

The four branches in Figure~\ref{fig:svk_plane_stress_lambda_z} show the through-thickness variation of the local out-of-plane adjustment. 
The two negative-$Y$ levels remain below unity, while the two positive-$Y$ levels remain above unity, consistent with Equation~\eqref{eq:lambda_z_SVK} and the local enforcement of $S_{33}=0$.
Each branch becomes nearly uniform over most of the beam, and the three axial discretisations approach the same interior levels.
Their visible differences are confined mainly to the short boundary-affected region, where the finer meshes provide a denser sampling of the local constitutive response.
For $\nu=0$, Equation~\eqref{eq:lambda_z_SVK} reduces to $\lambda_z=1$, so the non-trivial distribution is reported only for $\nu=0.49$.
Taken together, Figures~\ref{fig:svk_plane_stress_def}, \ref{fig:svk_plane_stress_lambda_t}, and \ref{fig:svk_plane_stress_lambda_z} show that coincident centrelines do not imply identical local cross-sectional deformation.

\FloatBarrier

\paragraph{Neo-Hookean model}
The same cantilever end-moment benchmark is considered using the compressible
Neo-Hookean plane-stress specialisation in \ref{app:nh_plane_stress}. The
calculations compare $\nu=0$ with the near-incompressible case
$\nu\approx0.49$. The global deformation is shown using
$n_{\mathrm{elem}}=20$, while the detailed local comparison uses
$n_{\mathrm{elem}}=100$.

\begin{figure}[htbp]
    \centering
    \includegraphics[width=0.8\textwidth]
    {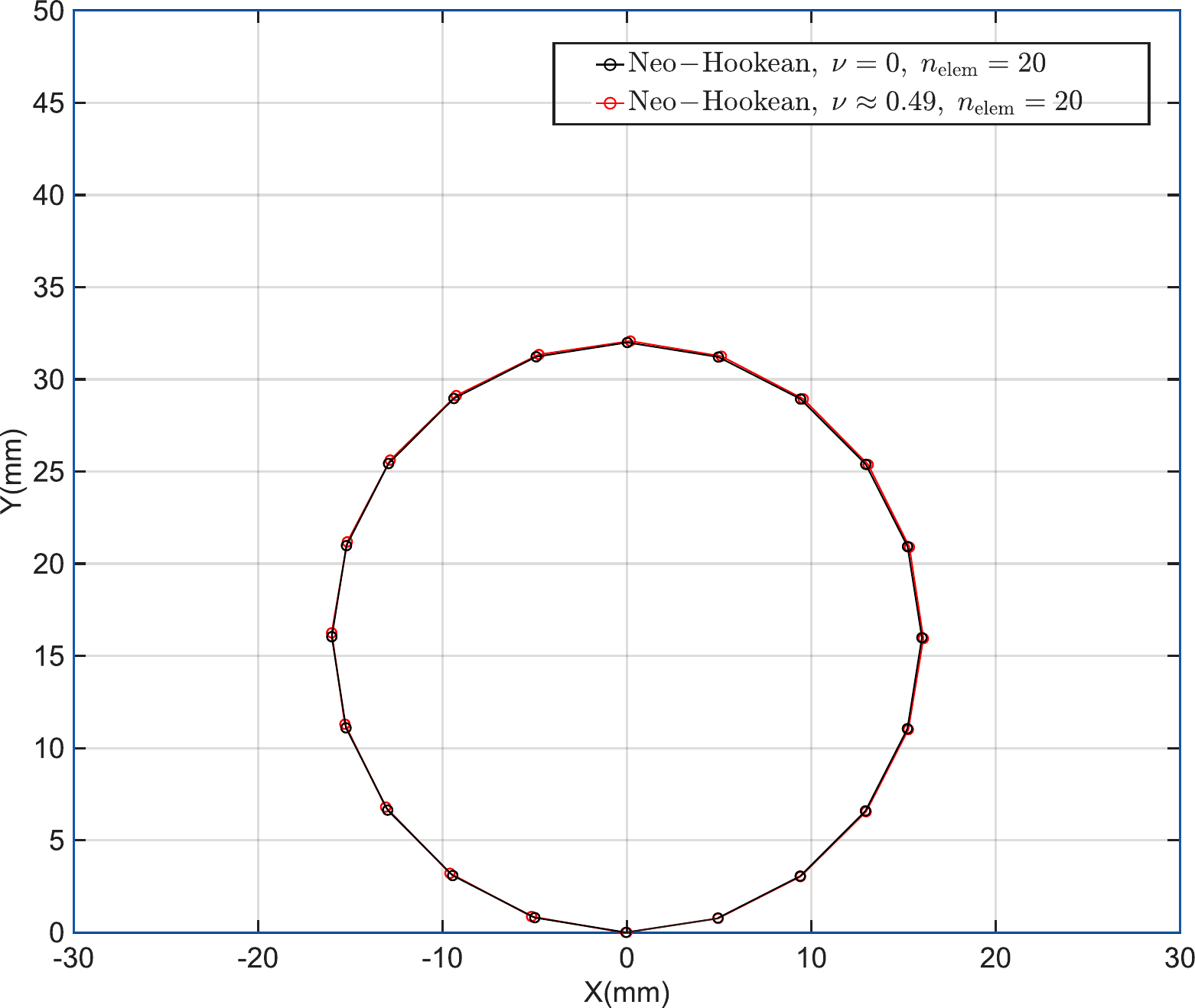}
    \caption{Deformed centreline for the Neo-Hookean plane-stress model with
    $\nu=0$ and $\nu\approx0.49$ using $n_{\mathrm{elem}}=20$. Open black
    and red circles denote $\nu=0$ and $\nu\approx0.49$, respectively.}
    \label{fig:nh_plane_stress_def}
\end{figure}

The two configurations in Figure~\ref{fig:nh_plane_stress_def} are practically indistinguishable over the complete beam.
The $20$-element discretisation is therefore sufficient for visualising the global configuration and the overlap of the two Poisson-ratio cases.
The Neo-Hookean centreline exhibits the same satisfactory mesh-refinement behaviour already established in the preceding Saint-Venant-Kirchhoff studies, so repeating the complete $20$-, $60$-, and $100$-element centreline comparison adds little information.

\begin{figure}[p]
    \centering

    \begin{subfigure}[t]{1\textwidth}
        \centering
        \includegraphics[width=\linewidth]
        {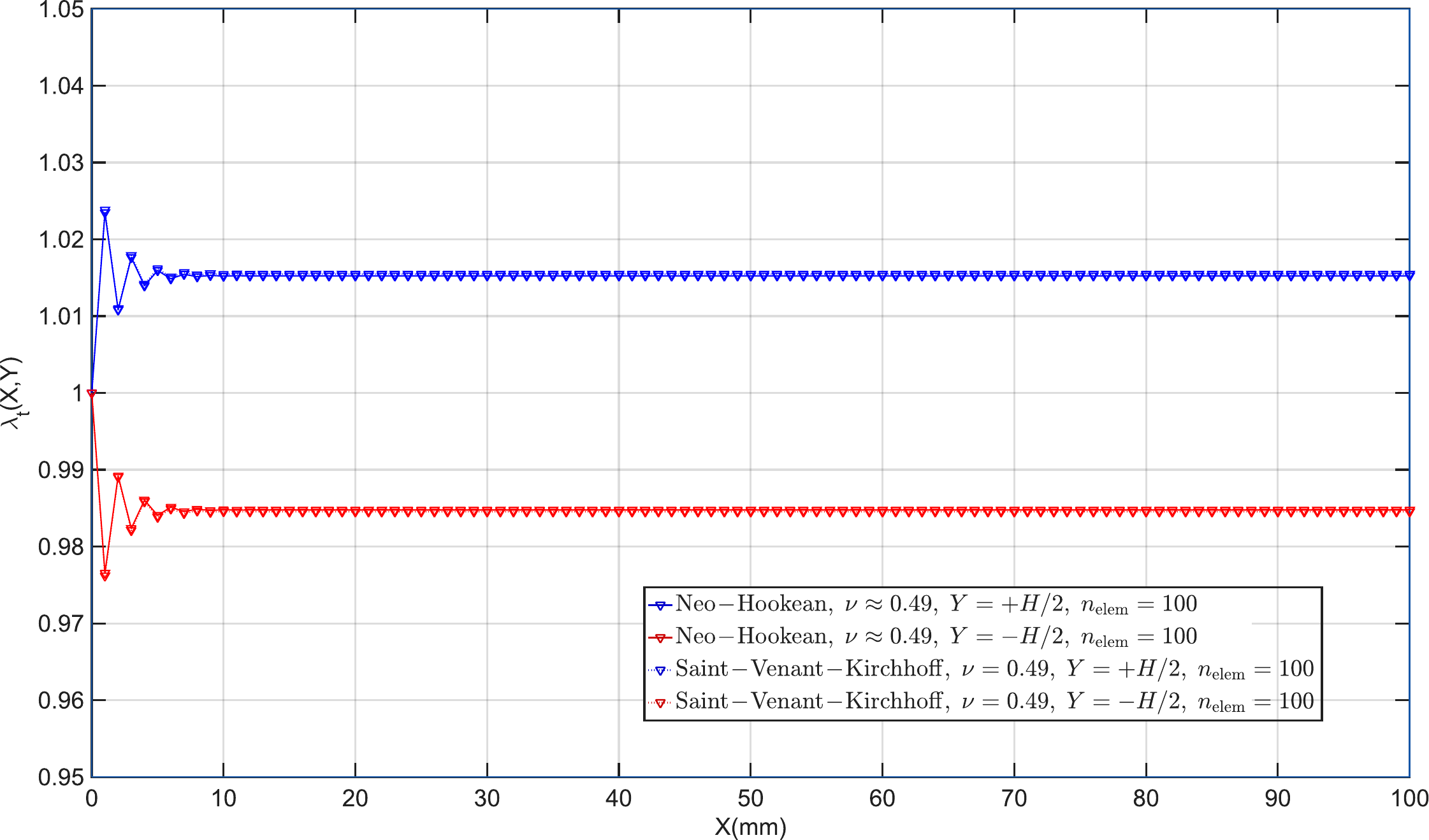}
        \caption{Surface values $\lambda_t(X,\pm H/2)$.}
        \label{fig:nh_plane_stress_lambda_t}
    \end{subfigure}

    \vspace{0.8em}

    \begin{subfigure}[t]{1\textwidth}
        \centering
        \includegraphics[width=\linewidth]
        {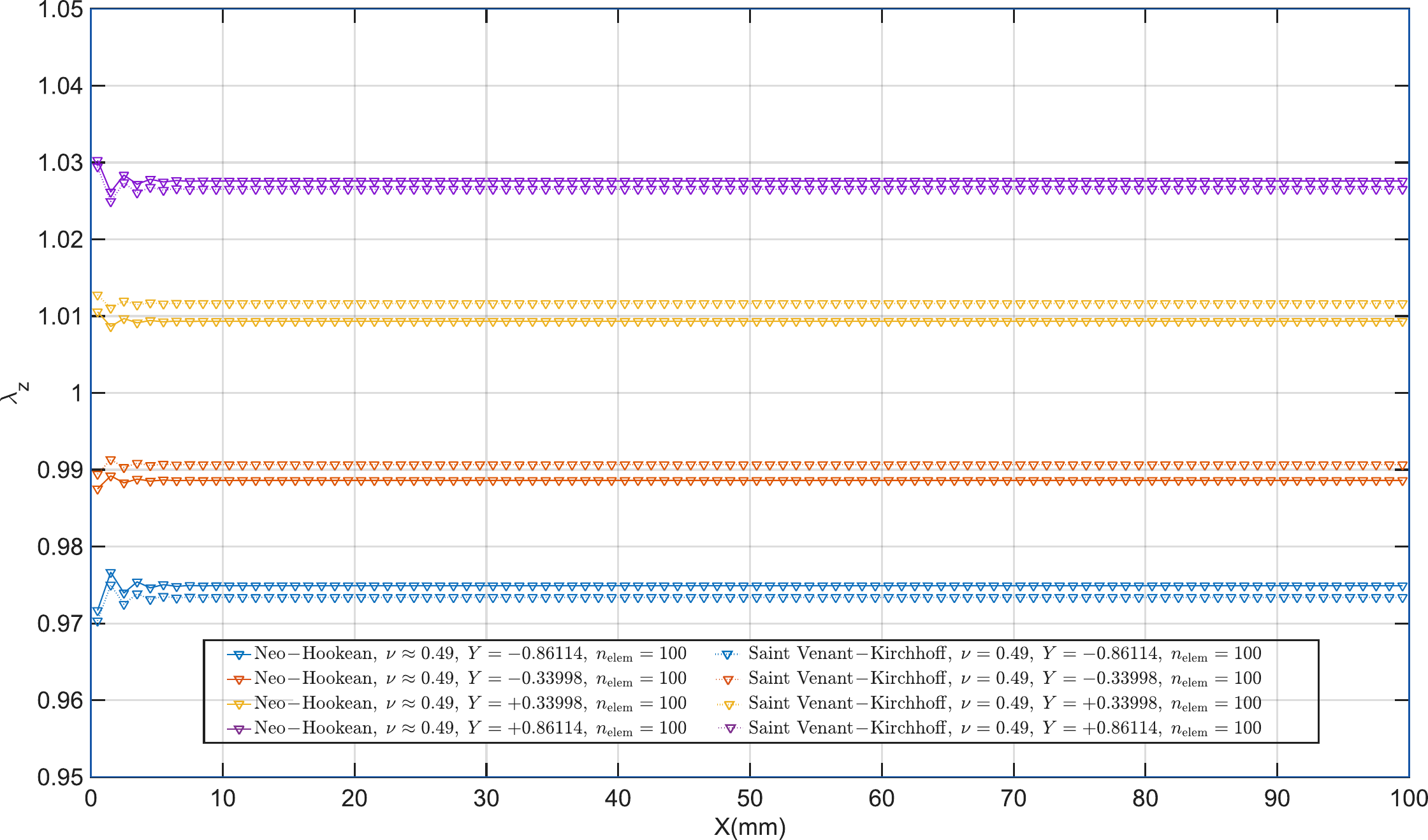}
        \caption{Out-of-plane stretch $\lambda_z$ at the four
        through-thickness integration-point levels.}
        \label{fig:nh_plane_stress_lambda_z}
    \end{subfigure}

    \caption{Comparison of the local plane-stress response of the
    Neo-Hookean and reduced Saint-Venant-Kirchhoff models using
    $n_{\mathrm{elem}}=100$:
    (a) surface values of the through-thickness scaling measure $\lambda_t$,
    and (b) out-of-plane stretch $\lambda_z$ at the four through-thickness
    integration-point levels. In (a), blue and red denote $Y=+H/2$ and
    $Y=-H/2$, respectively. In (b), blue, orange, yellow, and purple
    correspond to $Y=-0.86114$, $-0.33998$, $+0.33998$, and
    $+0.86114\,\mathrm{mm}$, respectively. Solid curves denote the
    Neo-Hookean model with $\nu\approx0.49$, while dotted curves denote the
    reduced Saint-Venant-Kirchhoff model with $\nu=0.49$. Open downward
    triangles identify the reported values.}
    \label{fig:nh_plane_stress_local}
\end{figure}

The local plane-stress response is compared directly with the corresponding reduced Saint-Venant-Kirchhoff result in Figure~\ref{fig:nh_plane_stress_local}.
The $100$-element discretisation is used for both constitutive models because it gives a cleaner representation of the local stretch fields and limits the visual influence of the short region near the clamped end.

Figure~\ref{fig:nh_plane_stress_lambda_t} shows essentially the same surface response for the two constitutive models.
The upper surface has $\lambda_t>1$ and the lower surface has $\lambda_t<1$.
A short oscillatory adjustment is visible near the clamped end, after which both branches settle rapidly to nearly uniform values.
The solid Neo-Hookean and dotted Saint-Venant-Kirchhoff curves are practically coincident over most of the beam, with only very small differences in the local clamped-end adjustment.
The use of $n_{\mathrm{elem}}=100$ is therefore not required by the global centreline response; it is used to resolve the local stretch field more cleanly and to reduce the prominence of the clamped-end effect.
The Neo-Hookean calculations show the same mesh-refinement trend as the preceding Saint-Venant-Kirchhoff results, so the full sequence of meshes is not repeated.

The comparison of $\lambda_z$ in Figure~\ref{fig:nh_plane_stress_lambda_z} resolves the constitutive response through the thickness more clearly. In both models, the two levels below the centreline give $\lambda_z<1$, whereas the two levels above the centreline give $\lambda_z>1$. After the short boundary adjustment, all four branches become almost uniform along the beam axis.

The Neo-Hookean and reduced Saint-Venant-Kirchhoff predictions remain close at all four integration-point levels, but small constitutive differences are visible in the settled values.
These differences are more apparent in $\lambda_z$ than in either the centreline configuration or the surface $\lambda_t$ response.
This is consistent with the distinct plane-stress reductions: the Neo-Hookean value follows from the nonlinear local equation in Equation~\eqref{eq:NH_plane_stress_equation}, whereas the reduced Saint-Venant-Kirchhoff relation is obtained analytically in \ref{app:svk_plane_stress}.
The observed separation is therefore local and constitutive rather than a difference in the global beam deformation.


\subsubsection*{Implications:}

The numerical comparisons show that the out-of-plane constitutive constraint can significantly influence the response as Poisson coupling becomes stronger.
In the plane-strain case, this constraint progressively affects the global response, whereas the plane-stress results in Figures~\ref{fig:svk_plane_stress_results}--\ref{fig:nh_plane_stress_lambda_z} permit transverse relaxation while leaving the global configuration essentially unchanged.
Related finite-strain beam models employ thickness-stretch variables under plane-strain assumptions \cite{HeLouDongKitipornchaiYang2018PlaneStrain,HeLouDongKitipornchaiYang2018Shearable}, while more general three-dimensional beam descriptions use extensible directors and enhanced strains to resolve cross-sectional deformation \cite{ChoiSauerKlinkel2021}.
The present approach instead partitions the two transverse responses within a planar setting: the in-plane thickness change is retained globally through $\gamma$, whereas the out-of-plane stretch $\lambda_z$ is governed solely by the selected constitutive constraint.
This allows plane strain and plane stress to be examined with identical global kinematics and, in the latter case, avoids an additional nodal field in the $Z$-direction.

The local enforcement of $S_{33}=0$ and the associated condensed tangent in Equations~\eqref{eq:plane_stress_condition} and \eqref{eq:L_condensed} follow the zero-stress assumption \cite{KlinkelGovindjee2002,AhmadiMcBrideSteinmannSaxena2025}, here coupled with the independently enriched in-plane beam kinematics.
The resulting formulation retains the transverse deformation relevant to the present planar setting while preserving a compact global finite-element description.

\section{Conclusions}
\label{sec:conclusions}

A geometrically-exact planar beam finite-element framework has been developed for large-deformation analysis with both linear and hyperelastic constitutive models.
While the conventional formulation with linear elastic material reproduces the three-dimensional reference deformation for the vertical-load benchmark and the analytical pure-bending solution for the end-moment benchmark, it is not directly applicable to hyperelastic beams requiring transverse deformation because its kinematics lack an independent thickness-stretch field and impose a plane-strain constraint.

The enriched formulation extends the beam kinematics by allowing independent deformation through the in-plane thickness while treating the out-of-plane direction through a separate constitutive reduction. However, the treatment of out-of-plane behaviour significantly affects the results obtained.
The plane-stress assumption determines the out-of-plane deformation by enforcing the zero-out-of-plane-stress constraint.
In contrast, the plane-strain assumption usually employed in nonlinear planar-beam finite element models suppresses out-of-plane deformation and introduces a reaction stress, producing a stiffer response for a different boundary-value problem. This difference reflects the imposed constraint rather than an intrinsically inaccurate plane-strain solution.
The end-moment calculations show the practical consequence of this distinction.
Under plane strain, the global curvature becomes increasingly sensitive to Poisson's ratio as the material approaches near-incompressibility.
Under plane stress, the global centrelines for the different Poisson's ratio cases considered are essentially coincident for both the Saint-Venant-Kirchhoff and Neo-Hookean models. This near-coincidence is consistent with the relatively small strains and mild constitutive nonlinearity encountered in the present examples, despite the large rotations.
Despite this similarity in the global deformation, the local in-plane and out-of-plane stretch responses remain non-trivial.
The two plane-stress material models predict very similar in-plane thickness changes, while small constitutive differences are more clearly resolved in the out-of-plane stretch through the beam thickness.

These results demonstrate that the cross-sectional response of a hyperelastic beam cannot be inferred from its centreline deformation alone.
Plane-strain and plane-stress conditions impose different three-dimensional constitutive restrictions and represent different boundary-value problems. Plane stress is appropriate for the present beam with an unconstrained out-of-plane direction, whereas plane strain represents a constrained two-dimensional structural idealisation. Experimental validation would therefore require different boundary conditions and experimental setups consistent with the respective assumptions.
The present framework provides a simple, yet consistent route for incorporating hyperelastic material models into geometrically exact beam finite elements. 
Future work will extend this framework to magnetomechanics, contact problems, and coupled fluid–structure interaction.

\appendix
\section{Material specialisations for the constitutive reductions}
\label{app:material_specialisations}

The plane-strain and plane-stress formulations in
Sections~\ref{sec:plane_strain} and~\ref{sec:plane_stress} are developed
independently of a particular hyperelastic strain-energy function. This
appendix collects the material-specific specialisations used in the numerical
examples of Section~\ref{sec:hyperelastic_results}. The
Saint-Venant-Kirchhoff model is considered under both plane strain and plane
stress, while the compressible Neo-Hookean model is considered under plane
stress to illustrate the local nonlinear determination of the out-of-plane
stretch.

\subsection{Saint-Venant-Kirchhoff material}
\label{app:svk_material}

The three-dimensional Saint-Venant-Kirchhoff (SVK) strain-energy density is
\begin{subequations}
\label{eq:SVK_material_relations}
\begin{equation}
\Psi_{\mathrm{SVK}}
=
\frac{\lambda_{\mathrm L}}{2}
\left(\tr\mathbf E\right)^2
+
\mu\,\mathbf E:\mathbf E,
\label{eq:SVK_energy}
\end{equation}
and the corresponding second Piola--Kirchhoff stress is
\begin{equation}
\mathbf S
=
\lambda_{\mathrm L}\tr(\mathbf E)\mathbf I
+
2\mu\mathbf E.
\label{eq:SVK_stress}
\end{equation}
\end{subequations}
The Lam\'e parameters are those defined in Equation~\eqref{eq:Lame_common}.

\subsubsection{Plane-strain specialisation}
\label{app:svk_plane_strain}

Under the plane-strain constraint of Section~\ref{sec:plane_strain},
$E_{33}=0$. The in-plane stress components therefore become
\begin{subequations}
\label{eq:SVK_plane_strain_stress_set}
\begin{equation}
\begin{aligned}
S_{11}
&=
(\lambda_{\mathrm L}+2\mu)E_{11}
+
\lambda_{\mathrm L}E_{22},
\\
S_{22}
&=
\lambda_{\mathrm L}E_{11}
+
(\lambda_{\mathrm L}+2\mu)E_{22},
\\
S_{12}
&=
2\mu E_{12},
\end{aligned}
\label{eq:SVK_plane_strain_inplane_stress}
\end{equation}
while the out-of-plane reaction stress is
\begin{equation}
S_{33}
=
\lambda_{\mathrm L}(E_{11}+E_{22}).
\label{eq:S33_SVK_plane_strain}
\end{equation}
\end{subequations}

Using the in-plane vectors defined in Equation~\eqref{eq:inplane_vectors}, the
plane-strain constitutive relation is
\begin{subequations}
\label{eq:SVK_plane_strain_constitutive_set}
\begin{equation}
\mathbf s
=
\bar{\mathbf C}_{\mathrm{SVK}}^{\mathrm{pe}}\mathbf e.
\label{eq:SVK_plane_strain_law}
\end{equation}
The reduced plane-strain matrix is
\begin{equation}
\bar{\mathbf C}_{\mathrm{SVK}}^{\mathrm{pe}}
=
\begin{bmatrix}
\lambda_{\mathrm L}+2\mu & \lambda_{\mathrm L} & 0\\
\lambda_{\mathrm L} & \lambda_{\mathrm L}+2\mu & 0\\
0 & 0 & \mu
\end{bmatrix},
\label{eq:SVK_plane_strain_matrix_Lame}
\end{equation}
or, equivalently,
\begin{equation}
\bar{\mathbf C}_{\mathrm{SVK}}^{\mathrm{pe}}
=
\frac{E_Y}{(1+\nu)(1-2\nu)}
\begin{bmatrix}
1-\nu & \nu & 0\\
\nu & 1-\nu & 0\\
0 & 0 & \dfrac{1-2\nu}{2}
\end{bmatrix}.
\label{eq:SVK_plane_strain_matrix}
\end{equation}
\end{subequations}
The corresponding reduced strain-energy density is
\begin{equation}
\Psi_{\mathrm{SVK}}^{\mathrm{pe}}
=
\frac12\mathbf e^{\mathrm T}
\bar{\mathbf C}_{\mathrm{SVK}}^{\mathrm{pe}}
\mathbf e.
\label{eq:SVK_plane_strain_energy}
\end{equation}
This specialisation makes the consequence of the plane-strain constraint
explicit. The in-plane response remains coupled through the
three-dimensional Lam\'e parameter, while $S_{33}$ acts as the reaction
required to maintain $E_{33}=0$.

\subsubsection{Plane-stress specialisation}
\label{app:svk_plane_stress}

For the SVK material, the plane-stress condition of
Section~\ref{sec:plane_stress} can be solved analytically. The out-of-plane
stress component is
\begin{subequations}
\label{eq:SVK_plane_stress_elimination}
\begin{equation}
S_{33}
=
\lambda_{\mathrm L}(E_{11}+E_{22}+E_{33})
+
2\mu E_{33}.
\label{eq:S33_SVK_plane_stress}
\end{equation}
Enforcing $S_{33}=0$ gives
\begin{equation}
E_{33}
=
-
\frac{\lambda_{\mathrm L}}{\lambda_{\mathrm L}+2\mu}
(E_{11}+E_{22})
=
-
\frac{\nu}{1-\nu}(E_{11}+E_{22}).
\label{eq:E33_SVK_plane_stress}
\end{equation}
Using $E_{33}=\tfrac12(\lambda_z^2-1)$, the corresponding positive
out-of-plane stretch is
\begin{equation}
\lambda_z
=
\sqrt{
1-
\frac{2\nu}{1-\nu}(E_{11}+E_{22})
}.
\label{eq:lambda_z_SVK}
\end{equation}
\end{subequations}
Thus, no local Newton iteration is required for this plane-stress
specialisation, provided the argument of the square root remains positive.

The in-plane second Piola--Kirchhoff stress components reduce to
\begin{subequations}
\label{eq:SVK_plane_stress_constitutive_set}
\begin{equation}
\begin{aligned}
S_{11}
&=
\frac{E_Y}{1-\nu^2}
(E_{11}+\nu E_{22}),
\\
S_{22}
&=
\frac{E_Y}{1-\nu^2}
(\nu E_{11}+E_{22}),
\\
S_{12}
&=
\frac{E_Y}{1+\nu}E_{12}.
\end{aligned}
\label{eq:SVK_plane_stress_inplane_stress}
\end{equation}
Using the in-plane vectors in Equation~\eqref{eq:inplane_vectors},
\begin{equation}
\mathbf s
=
\bar{\mathbf C}_{\mathrm{SVK}}^{\mathrm{ps}}\mathbf e.
\label{eq:SVK_plane_stress_law}
\end{equation}
The reduced plane-stress matrix is
\begin{equation}
\bar{\mathbf C}_{\mathrm{SVK}}^{\mathrm{ps}}
=
\frac{E_Y}{1-\nu^2}
\begin{bmatrix}
1 & \nu & 0\\
\nu & 1 & 0\\
0 & 0 & \dfrac{1-\nu}{2}
\end{bmatrix}.
\label{eq:SVK_plane_stress_matrix}
\end{equation}
\end{subequations}
The corresponding reduced strain-energy density is
\begin{equation}
\Psi_{\mathrm{SVK}}^{\mathrm{ps}}
=
\frac12\mathbf e^{\mathrm T}
\bar{\mathbf C}_{\mathrm{SVK}}^{\mathrm{ps}}
\mathbf e.
\label{eq:SVK_plane_stress_energy}
\end{equation}
Although Equation~\eqref{eq:SVK_plane_stress_matrix} has the familiar algebraic
form of the classical plane-stress elasticity matrix, it acts here on the
Green-Lagrange strain components and returns second Piola--Kirchhoff
stresses. The beam formulation therefore remains geometrically nonlinear.

\subsection{Neo-Hookean material under plane stress}
\label{app:nh_plane_stress}

The compressible Neo-Hookean model provides a nonlinear material
specialisation of the general plane-stress procedure in
Section~\ref{sec:plane_stress}. In this case, the out-of-plane stretch is
obtained from a nonlinear scalar equation at each material quadrature point.
The three-dimensional strain-energy density is
\begin{subequations}
\label{eq:NH_material_relations}
\begin{equation}
\Psi_{\mathrm{NH}}
=
\frac{\mu}{2}
\left(
\tr\mathbf C-3-2\ln J
\right)
+
\frac{\lambda_{\mathrm L}}{2}(\ln J)^2,
\label{eq:NH_energy}
\end{equation}
where
\begin{equation}
\mathbf C
=
\mathbf F^{\mathrm T}\mathbf F.
\label{eq:C_NH}
\end{equation}
The corresponding second Piola--Kirchhoff stress is
\begin{equation}
\mathbf S
=
\mu\left(\mathbf I-\mathbf C^{-1}\right)
+
\lambda_{\mathrm L}\ln J\,\mathbf C^{-1}.
\label{eq:NH_stress}
\end{equation}
\end{subequations}

For the block-diagonal plane-stress deformation gradient in
Equation~\eqref{eq:F_plane_stress},
\begin{subequations}
\label{eq:NH_out_of_plane_response}
\begin{equation}
C_{33}=\lambda_z^2,
\qquad
(C^{-1})_{33}=\frac{1}{\lambda_z^2},
\label{eq:NH_C33_relations}
\end{equation}
and hence
\begin{equation}
S_{33}
=
\mu\left(1-\frac{1}{\lambda_z^2}\right)
+
\lambda_{\mathrm L}
\frac{\ln(\lambda_zJ_{2D})}{\lambda_z^2}.
\label{eq:NH_S33}
\end{equation}
\end{subequations}
The plane-stress condition $S_{33}=0$ can therefore be written as
\begin{subequations}
\label{eq:NH_local_plane_stress_solve}
\begin{equation}
h(\lambda_z)
=
\mu(\lambda_z^2-1)
+
\lambda_{\mathrm L}\ln(\lambda_zJ_{2D})
=
0.
\label{eq:NH_plane_stress_equation}
\end{equation}
Its derivative is
\begin{equation}
h'(\lambda_z)
=
2\mu\lambda_z
+
\frac{\lambda_{\mathrm L}}{\lambda_z},
\label{eq:NH_hprime}
\end{equation}
which gives the local Newton update
\begin{equation}
\lambda_z^{(k+1)}
=
\lambda_z^{(k)}
-
\frac{
\mu\left[(\lambda_z^{(k)})^2-1\right]
+
\lambda_{\mathrm L}\ln\left(\lambda_z^{(k)}J_{2D}\right)
}{
2\mu\lambda_z^{(k)}
+
\dfrac{\lambda_{\mathrm L}}{\lambda_z^{(k)}}
}.
\label{eq:NH_lambda_Newton}
\end{equation}
\end{subequations}
A positive initial value is required so that the logarithm remains defined.
In an incremental calculation, the converged value from the preceding global
iteration or load step provides a natural initial estimate.

After convergence,
\begin{equation}
E_{33}
=
\frac12(\lambda_z^2-1),
\qquad
J
=
\lambda_zJ_{2D},
\label{eq:NH_E33_J}
\end{equation}
and the complete stress tensor is evaluated from Equation~\eqref{eq:NH_stress}.
The in-plane stress components are collected in the vector $\mathbf s$
defined in Equation~\eqref{eq:inplane_vectors}.

The material tangent with respect to the Green-Lagrange strain is
\begin{align}
L_{IJKL}
={}&
\lambda_{\mathrm L}
(C^{-1})_{IJ}(C^{-1})_{KL}
\nonumber\\
&+
\left(\mu-\lambda_{\mathrm L}\ln J\right)
\left[
(C^{-1})_{IK}(C^{-1})_{JL}
+
(C^{-1})_{IL}(C^{-1})_{JK}
\right].
\label{eq:NH_material_tangent}
\end{align}
The plane-stress tangent is obtained from the general condensation in
Equation~\eqref{eq:L_condensed}. Denoting the resulting matrix by
$\bar{\mathbf C}_{\mathrm{NH}}^{\mathrm{ps}}$, the local in-plane increment
has the form
\begin{subequations}
\label{eq:NH_plane_stress_element_set}
\begin{equation}
D\mathbf s
=
\bar{\mathbf C}_{\mathrm{NH}}^{\mathrm{ps}}D\mathbf e,
\label{eq:NH_condensed_matrix_law}
\end{equation}
and the material contribution to the element tangent is
\begin{equation}
\mathbf K_{\mathrm{mat,NH}}^{\mathrm{ps}}
=
\int_0^{L_e}\int_A
\mathbf B_E^{\mathrm T}
\bar{\mathbf C}_{\mathrm{NH}}^{\mathrm{ps}}
\mathbf B_E
\,\mathrm dA\,\mathrm dX.
\label{eq:NH_material_stiffness}
\end{equation}
\end{subequations}

For numerical implementation, each cross-sectional quadrature point is
treated locally. The fields $u$, $v$, $\beta$, and $\gamma$ and their axial
derivatives are first evaluated, followed by $\rho$, $\rho_{,X}$,
$\rho_{,Y}$, $\mathbf F_{2D}$, and $J_{2D}$.
Equation~\eqref{eq:NH_plane_stress_equation} is then solved for the positive
$\lambda_z$. The complete three-dimensional $\mathbf F$, $\mathbf C$, $J$,
$\mathbf S$, and $\mathbb L$ are evaluated at the converged state, after
which Equation~\eqref{eq:L_condensed} supplies the plane-stress tangent used in
the element residual and tangent matrix. The out-of-plane stretch remains an
internal constitutive quantity throughout and does not enlarge the global
finite-element system.

\section*{Acknowledgements}

\noindent The authors acknowledge support from the Engineering and Physical Sciences Research Council (EPSRC) under grant EP/Z535710/1. M. Hossain also acknowledges the Royal Society-NSFC International Exchange Grant IEC/NSFC/211316.

\bibliographystyle{elsarticle-num}
\bibliography{ref}

\end{document}